\documentclass[aos,seceqn,citesort,dvips]{arximspdf}
\usepackage{graphics}
\doi{10.1214/09-AOS779}
\volume{38}
\issue{4}
\pubyear{2010}
\firstpage{2242}
\lastpage{2281}

\makeatletter

\newcommand{\indep}{\perp\hspace*{-6.2pt}\perp}
\renewcommand{\bot}{\perp}

\newproclaim{example}{Example}[section]
\newproclaim{definition}{Definition}[section]
\newtheorem{lemma}{Lemma}[section]
\newtheorem{proposition}[lemma]{Proposition}
\newproclaim{condition}[lemma]{Condition}
\newproclaim{remark}{Remark}[section]
\newtheorem{theorem}[lemma]{Theorem}
\newtheorem{corollary}[lemma]{Corollary}
\newproclaim{algorithm}[lemma]{Algorithm}

\makeatother

\begin{document}
\begin{frontmatter}

\title{Stochastic kinetic models: Dynamic independence, modularity and
graphs\protect\thanksref{T1}}
\runtitle{Stochastic kinetic models}

\thankstext{T1}{Supported by the EPSRC and MRC (United Kingdom).}

\begin{aug}
\author[A]{\fnms{Clive G.} \snm{Bowsher}\corref{}\ead[label=e1]{C.Bowsher@statslab.cam.ac.uk}}
\runauthor{C. G. Bowsher}
\affiliation{University of Cambridge}
\address[A]{Centre for Mathematical Sciences\\
University of Cambridge\\
Wilberforce Road, Cambridge\\
United Kingdom\\
\printead{e1}} 
\end{aug}

\received{\smonth{8} \syear{2009}}
\revised{\smonth{12} \syear{2009}}

%
\begin{abstract}
The dynamic properties and independence structure of stochastic kinetic
models (SKMs) are analyzed. An SKM is a highly multivariate jump
process used to model chemical reaction networks, particularly those in
biochemical and cellular systems. We identify SKM subprocesses with the
corresponding counting processes and propose a directed, cyclic graph
(the kinetic independence graph or KIG) that encodes the local
independence structure of their conditional intensities. Given a
partition $[A,D,B]$ of the vertices, the graphical separation $A\perp
B|D$ in the undirected KIG has an intuitive chemical interpretation and
implies that $A$ is locally independent of $B$ given $A \cup D$. It is
proved that this separation also results in global independence of the
internal histories of $A$ and $B$ conditional on a history of the jumps
in $D$ which, under conditions we derive, corresponds to the internal
history of $D$. The results enable mathematical definition of a
modularization of an SKM using its implied dynamics. Graphical
decomposition methods are developed for the identification and
efficient computation of nested modularizations. Application to an SKM
of the red blood cell advances understanding of this biochemical
system.
\end{abstract}

%

\end{frontmatter}

\section{Introduction and summary}

The dynamic properties and conditional independence structure of
stochastic kinetic
models are analyzed using a marked point process framework. A stochastic
kinetic model or SKM is a highly multivariate jump process used to
describe chemical reaction networks. SKMs have become particularly
important as models of the network of interacting biomolecules in a
cellular system. The necessity of a stochastic process approach to the
dynamics of such
biochemical reaction systems is now clear
\cite{Wilkinson2009,Swain2008}, with SKMs
providing continuous-time, mechanistic descriptions firmly grounded in
chemical kinetic theory and the underlying statistical physics. The
Gillespie algorithm \cite{Gillespie76,Gillespie77} for simulation of
SKMs is now an important tool in the
science of systems biology. However, there are few analytical tools for
study of the dynamic properties of SKMs (although note
\cite{Kurtz2006,Goutsias2007}
and \cite{CGillespie2009}), especially when the SKM is
of modest or high dimension.

This paper develops what appear to be the first methods for analyzing
the local and
global dynamic independence structure implied by a given SKM and shows
how these may be used to uncover the modular architecture of the
network at coarser or finer levels of resolution. The required
information about the parameters of the SKM is
modest, and consistent with the partial information about these
currently available for many biochemical reaction networks. SKMs are
often thought of as continuous-time, homogeneous
Markov chains having nonfinite state space. However, the fact that
there are
a finite number of possible types of jump of the
process---corresponding to
the different types of possible biochemical reaction in the system---allows
formulation of both the SKM and its subprocesses as multivariate counting
processes. This turns out to be a fruitful approach for the problems
addressed here. In fact, the Markov property is not needed for the
results and methods of the paper. The main contributions may be
summarized as follows.

Graphical models for SKMs and dynamic molecular networks are
introduced. These kinetic independence graphs (KIGs) are directed,
cyclic graphs whose vertices are the different types
(or species) of biomolecule in the system. The KIG encodes local
independences that result from a lack of dependence of the conditional
intensity of a subprocess on the internal history of some of the species.

Given a partition $[A,B,D]$ of the vertices, the graphical separation
$A\perp B|D$ in the undirected version of the KIG has an intuitive
chemical interpretation and implies $A$ is locally (or
``instantaneously'') independent of $B$ given $A \cup D$ (and $B$ locally
independent of $A$ given $B \cup D$). It is proved that this separation
also results in conditional independence, over any finite time interval
$(0,t]$, of the internal histories of $A$ and $B$ conditional on a
history of the jumps in $D$. Conditions under which this history
corresponds to the internal history of $D$ are derived and are easily
checked computationally. Such a conditional independence is termed a
global (as opposed to local) dynamic independence here.

The new results enable mathematical definition of a modularization of
an SKM using its implied dynamics. Graphical decomposition methods are
developed for the identification of nested modularizations that allow
the extent of coarse-graining to be varied and provide computationally
efficient algorithms for large SKMs. Junction tree representations are
shown to provide a useful tool for visualizing, summarizing and
manipulating the modularizations. Applying the techniques of the paper
to an SKM that represents detailed empirical knowledge of the metabolic
network of the human red blood cell yields new insight into the
biological organization and dynamics of this cellular system.

Graphical models and their associated analytical and computational methods
allow the modularization of large, complex models into smaller components
and provide a particularly effective means of representing and analyzing
conditional independence relationships \cite{Lauritzen96,CDLS2007}.
Certain graphical approaches are now used quite extensively in
computational biology and have also
been readily assimilated by the wider biological scientific community,
which has long found
diagrammatic representations of reaction schemes useful
\cite{Huber2007}. However, rigorous graphical representations of biochemical
networks as dynamic processes---that is graphical models in the
statistical sense---do not appear to have been considered previously.

Indeed, graphical models for continuous time stochastic
processes in general are in an early stage of development. Didelez
\cite{Didelez2007,Didelez2008} introduced graphs based on the local independence
structure of
conditional intensities for finite state, composable Markov processes and
multivariate point processes, respectively; \cite{Nodelman2002} is an earlier
contribution, also for finite state Markov processes. SKMs require new methods
since interest is in dynamic independences between groups of species rather
than the counting processes for the different types of reaction per
se. Furthermore, the Markov process for species concentrations implied
by the SKM neither has finite state space, nor is it composable for
most SKMs of interest (see
Section \ref{sec3}).

In practice, the SKM is constructed from a large list of the
biochemical reactions that comprise the network under study. This list, or
``network reconstruction,'' is usually compiled using extensive experimental
evidence in the literature on the component parts of the system and their
molecular interactions \cite{Palsson2006}. Indeed, the approaches of
molecular biology and genetics, including genome sequencing, have
already proved
remarkably successful in providing life scientists with a very
extensive ``parts list'' for biology. Systems biology is an increasingly
influential, interdisciplinary approach that aims to describe
mathematically the stochastic dynamic behavior of the whole system as
an emergent property of the network of interacting biomolecules
\cite{Wilkinson2009}.

A principal challenge is thus to map from fine level descriptions such
as reaction lists and their implied SKMs to higher level,
coarse-grained descriptions of the dynamic properties. Related is the
increasingly held view that biochemical reaction networks are modular, that
is their architecture can be decomposed into units that perform ``nearly
independently'' \cite{Alon2005}, and that identifying such modules is a
crucial step in the endeavor to understand and, ultimately, to
selectively control
cellular systems. However, it is recognized that rigorous, mathematical
definition and identification of modularizations for biochemical
networks is
difficult, especially from a dynamic perspective (see
\cite{Stelling2006}, Chapter 3). As a result, such modularization
techniques have
been slow to develop, and there seems to be no prior work allowing for
stochastic and non-steady state dynamics. The dynamic independence
results and associated graphical methods developed here provide an
effective means of addressing these problems. Broadly speaking, the
paper also illustrates the utility of a statistical and probabilistic
approach to the dynamics of biological systems which, despite their
stochastic nature, have hitherto more often received the attention of
physical scientists.

The structure of the paper is as follows. Section \ref{sec21} introduces SKMs
and reaction networks in a manner requiring no previous background in
systems biology or biochemistry. Section \ref{sec22} defines an SKM as a marked
point process and provides a formal construction using the well-known
Gillespie algorithm as a point of departure. Section \ref{sec23} then shows how
to accommodate subprocesses of the SKM in a counting process framework
and discusses their conditional intensities and internal histories
(natural filtrations). Section \ref{sec3} introduces the kinetic independence
graphs, or KIGs, and examines local independence and graphical
separation in the undirected KIG. Section \ref{sec4} then relates these to
global conditional independence of species histories in Theorems
\ref{Main1} and \ref{Main2}, which are central to the paper. Rigorous
proofs of these theorems are quite involved and are given as Appendix
\ref{appA}. Section \ref{sec5} develops graphical decomposition methods and associated
theory for the identification of modularizations of SKMs, while Section
\ref{sec6} applies the techniques of the paper to the SKM of the human red blood
cell. Section \ref{sec7} highlights some directions for future research.

\section{SKMs and counting processes}\label{sec2}
\subsection{Introducing the SKM and reaction networks}\label{sec21}

A stochastic kinetic model is a continuous-time jump process modeling the
state of a chemical system, $X(t)=[X_{1}(t),\ldots,
X_{n}(t)]^{\prime}$, where $%
X_{i}(t)$ is interpreted as the nonnegative, integer number of
molecules of
type $i$ present at time $t$. The set of different types of molecule or the
\textit{species set} is given by $\mathcal{V}:=\{1,\ldots,n\}$. There are a
finite number of possible types of jump in $X(t)$ that may take place,
corresponding to the different types of possible \textit{reaction},
$m\in
\mathcal{M}:=\{1,\ldots,M\}$. It is particularly useful for our purposes to
view an SKM as a marked point process or MPP in which the points or
``events''
correspond to the jump times of the process $X(t)$. Mathematically, a
particular reaction can then be identified with an element of the finite
mark space and each mark indicates the type of jump associated with the
corresponding
jump time.

An SKM is denoted here by $\{T_{s},Z_{s}\}_{s\geq1}$, where $%
T_{s}$ is the $s$th jump time. The mark $Z_{s}\in\{S_{m}|m\in\mathcal
{M}\}$
is the value of the jump and is interpreted as the \textit{changes} in the
number of molecules of each species. The matrix $S:=[S_{1},S_{2},\ldots,S_{M}]$
is usually known as the stoichiometric matrix. Any two columns of $S$ are
taken to be nonequal; hence, there is a bijection between the mark
space and
$\mathcal{M}$. A formal construction of an SKM is given below in
Section \ref{sec22} but it is
helpful at this stage to note the following linear equation determining the
dynamic evolution of $X(t)$:
\[
X(t)=X(0)+SN(t),\qquad t\geq0,
\]
where $N(t)=[N_{1}(t),\ldots,N_{M}(t)]^{\prime}$ is the $M$-variate counting
process associated with the marked point process $\{T_{s},Z_{s}\}
_{s\geq1}$.
Thus, $N_{m}(t)$ is interpreted\vspace*{1pt} as counting the number of reactions of
type $m$
during $(0,t]$. Denote by $\mathcal{F}_{t}^{N}:=\sigma(N(s)$; $0\leq
s\leq t)$ the internal
history of the entire process and by $\mathcal{F}_{t}^{m}:=\sigma
(N_{m}(s)$; $0\leq s\leq t)$ the internal history of the $m$th counting
process. The probability law of $%
N(t)$, and hence that of $X(t)$, is determined by what are known as the
$%
\mathcal{F}_{t}^{N}$-conditional intensities, $[\lambda_{m}(t);m\in
\mathcal{%
M}]$.

The conditional intensity concept is important for an understanding of the
paper. At time $t$, each intensity $\lambda_{m}(t)$ is interpreted as
the \textit{local}
(or \textit{instantaneous}) \textit{rate} of reaction $m$, conditional on the
internal history of \textit{the entire process }$\mathcal{F}_{t}^{N}$.
Confining attention to a finite interval of time $\mathcal{T}$,
provided that
$N(t)$ has finite expectation $\forall t\in\mathcal{T}$ (and that
$[\lambda
_{m}(t);t\in\mathcal{T}]$ is bounded by an integrable random
variable), each intensity
is a local rate of reaction in exactly the chemical sense---that is, $%
\lambda_{m}(t+)=\lim_{h\downarrow0}\mathsf{E}[h^{-1}\{
N_{m}(t+h)-N_{m}(t)%
\}|\mathcal{F}_{t}^{N}]$, the conditionally expected number of
reactions of type
$m$ per unit time in the limit as $h$ goes
to zero. Of course, the intensities are themselves random variables
(r.v.'s) since the evolution
of $N$ up to time $t$ is itself a stochastic process, hence the appearance
of the conditional expectation. A technical subtlety is that $\lambda_{m}(t)$
is defined to have sample paths that are left-continuous (with limits from
the right), compared to the right-continuous sample paths of $X(t)$. A
heuristic chemical interpretation is that if a jump in $X$ takes place
at $t$, then future jumps are (locally) determined by the intensity evaluated
``immediately after'' $t$.

A basic familiarity with the chemical representation and interpretation of
reactions is helpful in what follows (see also \cite{Wilkinson2009} for
an accessible introduction). Each
reaction $m\in\mathcal{M}$ has the chemical representation
%
%
\begin{equation} \label{reaction}
\sum_{i\in R[m]}\alpha_{i}X_{i}\rightarrow^{m}\sum_{j\in P[m]}\beta
_{j}X_{j},
\end{equation}
which is read as follows: when reaction $m$ takes place, $\alpha_{i}$
molecules of type $i$ are consumed for each $i$ in the subset
$R[m]\subset
\mathcal{V}$, and $\beta_{j}$ molecules of type $j$ are produced for
each $j
$ in the subset $P[m]\subset\mathcal{V}$. The species $R[m]$ are
called the
\textit{reactants} (or inputs) of the reaction $m$, and the species $P[m]$
are called the \textit{products} (or outputs) of $m$. The integer
coefficients $[\{\alpha_{i}\},\{\beta_{j}\}]$ are known as the
stoichiometries of the reaction. If a species $k$ is a reactant but not a
product, then its corresponding entry in the stoichiometric matrix $S$
(i.e., the change
in the level of $k$ caused by reaction $m$) is
given by $S_{km}=-\alpha_{k}$. Alternatively, if species $k$ is a product
but not a reactant, then $S_{km}=\beta_{k}$. There is no assumption
that $%
R[m]\cap P[m]=\varnothing$, and if $k$ is both a product and a reactant
then $%
S_{km}=\beta_{k}-\alpha_{k}$. A common situation in this case is $%
\beta_{k}=\alpha_{k}$, that is $k$ acts as a ``catalyst,'' increasing the
rate of the reaction but not itself being ``changed'' by the
reaction---that is, not
itself being overall consumed or produced when $m$ takes place.
Formally, the sets $R[m]
$ and $P[m]$ are defined by allowing zero stoichiometries and writing
the $m$th reaction as $\sum_{i\in\mathcal{V}}\alpha_{i}X_{i}\rightarrow
^{m}\sum_{j\in\mathcal{V}}\beta_{j}X_{j}$. Then $R[m]:=\{i\in\mathcal
{V}%
|\alpha_{i}>0\}$ and $P[m]:=\{j\in\mathcal{V}|\beta_{j}>0\}$.

In systems biology, a living cell is often viewed as a network of interacting
biomolecules of different types, with $n$ and $M$ both large (and often
$M>n$%
). The interaction is selective---only species that are reactants for some
reaction $m$ can together react to give products. Each reaction involves
only a few species, so the cardinality of $R[m]\cup P[m]$ is small. Certain
reactions are ``coupled'' in that a product of one reaction is also a reactant
of another reaction. From a stochastic process perspective, the specification
of the list of component reactions as in (\ref{reaction})
for all $m \in\mathcal{M}$ implies dependences between the levels (or
concentrations) of
the different biomolecules.

As a simple but nonetheless biochemically meaningful illustration, consider
the following example of an SKM.
\begin{example}
\label{Example1_WellBeh}
Consider the SKM with the 5 different species $\mathcal{V}%
=\{P$, $R$, $g$, $P_{2}$, $gP_{2}\}$ and the 6 reactions
\begin{eqnarray*}
g&\rightarrow^{\mathrm{trc}}&g+R,\qquad R\rightarrow^{\mathrm{trl}}R+P,\qquad 2P\rightarrow
^{d}P_{2},
\\
P_{2}&\rightarrow^{rd}&2P,\qquad g+P_{2}\rightarrow^{b}gP_{2},\qquad
gP_{2}\rightarrow^{ub}g+P_{2}.
\end{eqnarray*}
\end{example}

The gene $(g)$ is responsible for the production of molecules of
protein $(P)$ via the intermediate (mRNA) species $(R)$. In this simplified
representation, $g$ and $R$ act as simple catalysts in the reactions $\mathrm{trc}$
(``transcription'') and $\mathrm{trl}$ (``translation''), respectively. The third
reaction $d$ consists of the binding of 2 molecules of $P$ (the sole
reactant) to form the
new molecule $P_{2}$ (the sole product). The fourth reaction $rd$ is
the reverse of the third.
The fifth reaction sets up a ``negative feedback cycle'' whereby the
production of
$P$ is negatively self-regulated by the binding of $P_{2}$ to $g$ to form
the distinct species $gP_{2}$. Genes bound in this way to $P_{2}$ are not
then available to participate in the $\mathrm{trc}$ reaction, thus preventing
over-production of the protein. We shall return later to the same example.

\subsection{Defining and constructing the SKM}\label{sec22}

The Gillespie stochastic simulation algorithm
\cite{Gillespie76,Gillespie77} has become an important tool in
biological science for
studying biochemical and cellular systems. Given its familiarity in
mathematical and computational biology, the following construction of
an SKM
as a marked point process takes as its point of departure the conditional
distributions employed in the Gillespie algorithm. For our purposes,
the algorithm is
usefully viewed as outputting a realization of the MPP $\{T_{s},Z_{s}\}$,
from which the resultant process $X(t)$ is easily constructed as in
(\ref{Xt}) below. Readers less concerned with formal constructions and
already familiar with stochastic kinetics may proceed safely to
Section~\ref{sec23} after noting Definition \ref{SKM} of an SKM and (\ref{cmgm}) for
the conditional reaction intensities (or ``hazards'').

Denote the numbers of molecules of all species at time $T_{s}$ by
$Z_{s}^{X}:=Z_{s-1}^{X}+Z_{s}$ $(s=1,2,\ldots)$. Let $Z_{0}^{X}$ be the
initial, deterministic state of the system, and define $T_{0}:=0$. We
write the $\sigma$-field generated by the first $r$ points and marks as
$\mathcal{F}_{T_{s}}:=\sigma(T_{r},Z_{r};r=1,\ldots,s)$. Also let $\mathcal
{F}_{T_{s+1^{-}}}:=\sigma
(T_{s+1},T_{r},Z_{r};r=1,\ldots,s)$, where the $(s+1)$th mark is excluded
from the generating collection of random variables.

Now introduce the important propensity (or reaction rate) function for
the $m$th reaction, $\lambda_{m}(Z_{s}^{X})$,
where $\lambda_{m}\dvtx\mathbb{N}_{0}^{n}\rightarrow\lbrack0,\infty)$ is
continuous. The conditional distributions implied by stochastic
kinetic theory
\cite{Gillespie92} and employed in the Gillespie algorithm are given by
%
%
\begin{eqnarray}\label{Wait}
\mathsf{P}(T_{s+1}>t|\mathcal{F}_{T_{s}})=\exp
\Biggl\{-(t-T_{s})\sum_{m=1}^{M}\lambda_{m}(Z_{s}^{X})\Biggr\},\nonumber\\[-8pt]\\[-8pt]
\eqntext{t>T_{s},s=0,1,\ldots,}
\end{eqnarray}
that is the waiting time to the next occurrence of a jump (reaction) is
exponentially distributed
with parameter $\sum_{m=1}^{M}\lambda_{m}(Z_{s}^{X})$; and
%
%
\begin{equation} \label{Jump}
\mathsf{P}(Z_{s+1}=S_{m}|\mathcal{F}_{T_{s+1^{-}}})=\lambda
_{m}(Z_{s}^{X})\Big/\sum_{m=1}^{M}\lambda_{m}(Z_{s}^{X}),
\end{equation}
which gives the mark (or jump) distribution.
Note that both the waiting time and mark distributions depend only
on $Z_{s}^{X}$, the levels of the species present following the $s$th
reaction. The pure jump process $X(t)$ is given straightforwardly, for
$t\geq0$, by
%
%
\begin{equation} \label{Xt}
X(t):=Z_{\max\{s\dvtx T_{s}\leq t\}}^{X},\qquad X(0):=Z_{0}^{X},
\end{equation}
it being well known that $X(t)$ is a time-homogeneous Markov chain
under $%
\mathsf{P}$.

It turns out to offer significant advantages and simplification to
adopt a
MPP framework for the problems addressed in the paper. An SKM is thus
defined here directly in terms of the MPP $\{T_{s},Z_{s}\}$ and its
corresponding counting processes. It is implicit in our definition of a MPP
that $T_{s}<T_{s+1}$ whenever $T_{s}<\infty$ $(s\geq1)$. Thus,
reactions occur
instantaneously and no two reactions ever have identical occurrence
times in continuous time.
The physical interpretation is that reaction durations are negligible
and may be ignored. The random variables $T_{s}$ are $(0,\infty
]$-valued, with the interpretation that
less than $s$ reactions take place during the time interval $[0,\infty
)$ if $T_{s}=\infty$.
The flexibility gained will not be needed routinely, but may be useful
for cellular systems that can enter an inactive or quiescent state. The
stability
condition $\lim_{s\rightarrow\infty}T_{s}=\infty$ a.s. is
imposed, which
is equivalent to the statement that only finitely many reactions occur in
any finite time interval (sometimes known as nonexplosivity).
\begin{definition}
\label{SKM}
A \textit{stochastic kinetic model} (\textit{SKM}) is the
MPP\break $[\{T_{s},Z_{s}\}_{s\geq1}$, $S,\mathsf{P}]$ with mark space given by the
columns of $S$, $\{S_{m}|m\in\mathcal{M}\}$, where no 2 columns of $S$ are
equal; and where the probability measure $\mathsf{P}$ is such that
(\ref{Wait}) holds $\mathsf{P}$-a.s.
on $\{T_{s}<\infty\}$, (\ref{Jump})
holds $\mathsf{P}$-a.s. on $\{T_{s+1}<\infty\}$, and $\lim_{s\rightarrow
\infty}T_{s}=\infty$ $\mathsf{P}$-a.s.

Equivalently, the SKM may be denoted by the corresponding multivariate
counting process (MVCP), $[N,S,\mathsf{P}]$, where $N:=[N_{m}(t);m\in
\mathcal{M}]_{t\geq0}$, and $N_{m}(t)=\sum_{s\geq1}1(T_{s}\leq
t)1(Z_{s}=S_{m})$ counts the number of reactions of type $m$ that occur
during $[0,t]$.
\end{definition}

Note that by definition the reaction counting processes $%
\{N_{m}(t);m\in\mathcal{M}\}$ have no jump times in common.
If the stability condition $T_{s}\rightarrow\infty$ a.s. holds,
there exists for any propensity functions ${\lambda_{m}(X_{t-})}$---see
\cite{Jacobsen82}, Theorem 1.7, page 56---a unique or \textit
{canonical} SKM satisfying Definition \ref{SKM} on $(\Omega,\mathcal
{F})$, where $\Omega$ is the space of $M$-variate counting process
paths (\cite{Jacobsen82}, Definition 1.2, page 53), $N$ is the
identity map
from $\Omega\rightarrow\Omega$, and $\mathcal{F}=\sigma(N(t);t\geq0)$.

It follows from (\ref{Wait}) and (\ref{Jump}) that the
propensity functions give the $%
\mathcal{F}_{t}^{N}$-conditional intensity process $\lambda(t)$ in the MVCP
sense (see \cite{Jacobsen82}, Definition 2.7), that is, $%
\lambda(t)=[\lambda_{m}(X_{t-})]_{m\in\mathcal{M}}$. When $N(t)$ has
finite expectation
$\forall t>0$, this means that $[N_{m}(t)-\int_{0}^{t}\lambda
_{m}(s)\,ds]$ is an
$\mathcal{F}_{t}^{N}$-martingale $\forall m$. That the intensities satisfy
%
%
\begin{eqnarray} \label{localEv}
\lim_{h\downarrow0}\frac{1}{h}\mathsf{P}\bigl(N_{m}(t+h)-N_{m}(t)
&=&1|\mathcal{F%
}_{t}^{N}\bigr)=\lambda_{m}(X_{t}),\qquad m\in\mathcal{M},\nonumber\\[-8pt]\\[-8pt]
\lim_{h\downarrow0}\frac{1}{h}\mathsf{P}\bigl(\bar{N}(t+h)-\bar{N}(t) &>&1|%
\mathcal{F}_{t}^{N}\bigr)=0, \nonumber
\end{eqnarray}
where $\bar{N}(t):=\sum_{m\in\mathcal{M}}N_{m}(t)$, is in fact a
principle conclusion of the arguments of stochastic kinetic theory
\cite{Gillespie92}. The assumptions of the theory are that the system is
spatially homogeneous (or ``well-stirred''), confined to a fixed volume and
held at constant temperature. Under these assumptions, (\ref{Wait}) and
(\ref{Jump}) have a firm physico-chemical basis \cite{Gillespie92}.

It plays a significant role in what follows that the
theory implies that the $\mathcal{F}_{t}^{N}$-intensities, $\lambda_{m}(t)$,
have the form
%
%
\begin{equation}
\label{cmgm}
\lambda_{m}(t)=c_{m}g_{m}\bigl\{X^{R[m]}(t-)\bigr\},
\end{equation}
where $c_{m}>0$ is a deterministic (``rate'') constant, and $g_{m}\{\cdot
\}\geq0$ is a
continuous function depending \textit{only on the levels of the
reactants} $R[m]$.


\subsection{SKM subprocesses---histories and intensities}
\label{sec23}

For any subset of\break molecular species $A\subseteq\mathcal{V}$, let the vector
process $\{X^{A}(t)\}:=\{X_{i}(t);i\in A\}$ denote the corresponding
subprocess of $X$. We identify $X^{A}$ with its MVCP, analogously to the
treatment of $X=X^{\mathcal{V}}$ above. For
$A\subseteq\mathcal{V}$, consider the subset of reactions $\Delta
(A)\subseteq\mathcal{M}$ that change (the level of) $A$, that is,
$\Delta(A):=\{m\in
\mathcal{M}\dvtx S_{m}^{A}\neq\mathbf{0}\}$, where $S_{m}^{A}$ is the subvector
of $S_{m}$ corresponding to the elements of $A$. One can identify $X^{A}$
with the MPP, $\{T_{s}^{A},Z_{s}^{A}\}$, where each jump time $T_{s}^{A}$
corresponds to the occurrence of some reaction in $\Delta(A)$; the
mark $%
Z_{s}^{A}$ gives the resultant jumps in the elements of $A$ and takes
its value in the mark space $E^{A}:=\{S_{m}^{A}|m\in\Delta
(A)\}$. This results in the following definition of an SKM subprocess.
\begin{definition}
\label{NADeltaA} Let $[N,S,\mathsf{P}]$
be an SKM and for $A\subset\mathcal{V}$, let $\Delta(A)$ be the
nonempty, finite subset $\{m\in\mathcal{M}\dvtx S_{m}^{A}\neq\mathbf{0}\}$.
Denote by $\mathcal{M}(\Delta(A))$ the partition of $\Delta(A)$ obtained
by grouping reactions that change $A$ identically, that is by applying
to $%
\Delta(A)$ the equivalence relation
\[
m\sim_{A}m^{\prime}\quad\Leftrightarrow \quad S_{m}^{A}=S_{m^{\prime}}^{A}.
\]
Denote the $e$th element of $\mathcal{M}(\Delta(A))$ by $\mathcal{M}%
_{e}(\Delta(A))$, $e=1,\ldots,|\mathcal{M}(\Delta(A))|$. The
\textit{subprocess of the SKM},
$N^{A}(t)$, is the $|\mathcal{M}(\Delta(A))|$-variate counting
process given by
%
%
\begin{equation} \label{NADeltaAEq}
N^{A}(t):= \biggl\{ \sum_{m\in\mathcal{M}_{e}(\Delta(A))}N_{m}(t)
\biggr\}
_{e=1,\ldots,|\mathcal{M}(\Delta(A))|}.
\end{equation}
The internal history of $N^{A}(t)$ is denoted $\mathcal{F}_{t}^{A}$.
\end{definition}

Note that since $\mathcal{M}(\Delta(A))$ is a partition of $%
\Delta(A)$, the components of $N^{A}(t)$ have no jumps in common. Each
element of the MVCP $N^{A}(t)$ thus counts the number of times
reactions in $\Delta(A)$ have occurred that result in a given change
in $A$. Intuitively, putting these elements together for all possible
types of change in $A$ to form a sample path of $N^{A}(t)$ captures
exactly the ``information'' given by the corresponding sample path of
$X^{A}(t)$. Indeed, there is
a bijection between the sample paths of $N^{A}(t)$ and those of
$X^{A}(t)$. The following technical lemma establishes that the internal
history of the MVCP $%
N^{A}(t)$ is identical to that of $X^{A}(t)$.
\begin{lemma}
\label{Filts}For $A\subseteq\mathcal{V}$, let $N^{A}(t)$ be a subprocess
of an SKM as in Definition \ref{NADeltaA} and let $\mathcal{F}%
_{t}^{X_{A}}:=\sigma(X^{A}(s);s\leq t)$ be the internal history of the jump
process $X^{A}(t)$. Then $\mathcal{F}_{t}^{A}=\mathcal{F}_{t}^{X^{A}}$ $
\forall t\geq0$. Furthermore, if $[A,B,D]$ is a partition of $\mathcal{V}$
then $\mathcal{F}_{t}^{N}=\mathcal{F}_{t}^{A}\vee\mathcal
{F}_{t}^{B}\vee
\mathcal{F}_{t}^{D}=\mathcal{F}_{t}^{X}$, $\forall t\geq0$.
\end{lemma}
\begin{pf}
A proof that $\mathcal{F}_{t}^{A}=\mathcal{F}_{t}^{X^{A}}$ is given in
Appendix \ref{Prooffs}. For $\mathcal{F}_{t}^{N}=\mathcal{F}_{t}^{X}$,
take $%
A=\mathcal{V}$. Finally, $\mathcal{F}_{t}^{X}=\mathcal
{F}_{t}^{X^{A}}\vee
\mathcal{F}_{t}^{X^{B}}\vee\mathcal{F}_{t}^{X^{D}}$ since $%
X(t)=[X^{A}(t)^{\prime},X^{B}(t)^{\prime}$, $X^{D}(t)^{\prime
}]^{\prime}$.
\end{pf}

One advantage of a counting process definition of the subprocess for the
species in $A\subset\mathcal{V}$ is that one may speak of the $\mathcal
{F}%
_{t}^{N}$-intensity for the subprocess and interpret this in the usual
manner as determining the \textit{local} or \textit{instantaneous}
dependence of the subprocess on the full internal history of the SKM, $%
\mathcal{F}_{t}^{N}$.
\begin{proposition}
\label{Ints}For $A\subseteq\mathcal{V}$, let $N^{A}(t)$ be a
subprocess of
an SKM as in Definition \ref{NADeltaA}. The
$\mathcal{F}_{t}^{N}$-conditional intensity
under $\mathsf{P}$ is given by
%
%
\begin{equation} \label{SupInt}
\lambda^{A}(t):= \biggl\{ \sum_{m\in\mathcal{M}_{e}(\Delta(A))}\lambda
_{m}(t) \biggr\}_{e=1,\ldots,|\mathcal{M}(\Delta(A))|}.
\end{equation}
\end{proposition}
\begin{pf}
Immediate from (\ref{NADeltaAEq}) on noting that the
intensities of the
superpositions of the counting processes are the sums of the corresponding
intensities.
\end{pf}

\noindent Notice that each element of the intensity, $\lambda
_{e}^{A}(t)$, is
the sum of the intensities (or stochastic rates) of all those reactions
that result in the corresponding change in $A$.

It follows from the equations in (\ref{localEv}) that, for any $%
A\subseteq\mathcal{V}$, the probability conditional on $\mathcal{F}_{t}^{N}$
that, during $(t,t+h]$, there is no change in $X^{A}$ is equal to $%
1-h\sum_{e=1}^{|\mathcal{M}(\Delta(A))|}\lambda_{e}^{A}(t)+o(h)$.
Similarly, the probability conditional on $\mathcal{F}_{t}^{N}$ that, during
$(t,t+h]$, there is exactly one jump in $X^{A}$ equal to $S_{m}^{A}$, for
some $m\in\mathcal{M}_{e}(\Delta(A))$, and also that no other
reaction $%
m^{\prime}\in\mathcal{M}$ occurs is equal to $h\lambda_{e}^{A}(t)+o(h)$
for $e=1,\ldots,|\mathcal{M}(\Delta(A))|$. Summing over all of the foregoing,
mutually exclusive events shows that these have conditional probability
equal to $1+o(h)$. Thus, in this infinitesimal sense, the $\mathcal{F}%
_{t}^{N}$-intensity $\lambda^{A}(t)$ may be interpreted as determining the
local dependence of $N^{A}(t)$ on $\mathcal{F}_{t}^{N}$.

\section{Kinetic independence graphs}
\label{sec3}

The identification of subprocesses of the SKM with their corresponding
MVCPs (see Section \ref{sec23}) greatly facilitates the construction of a
\textit{kinetic independence graph} encoding the local independence
structure of the SKM---see Definition \ref{LIGnew} below. The use of
the local
independence concept in constructing graphical models for continuous time
processes owes much to Didelez \cite{Didelez2007,Didelez2008}.
However, SKMs require new methods since interest is in dynamic independences
between groups of species rather than the reaction counting processes $%
[N_{m}(t)]_{m\in\mathcal{M}}$ per se. Thus, the vertex set of the
graph will be $\mathcal{V}$ rather than $\mathcal{M}$.

It is worth noting that existing graphical models for continuous-time Markov
chains \cite{Nodelman2002,Didelez2007} are not applicable to SKMs
because the Markov process $X(t)$ neither has finite state space, nor
is it
composable for most SKMs of interest. Roughly speaking, composability
\cite{Didelez2007} implies that any change of state in $X(t)$ can be
represented as a change in only one of several components.
Consider the use of $X^{A}(t)$ and $X^{\mathcal{V}%
\setminus A}(t)$ as components [since if $X(t)$ is composable with more
than 2 subsets of species as components, it must be composable with
just 2
components]---either the paths of $X^{A}(t)$ and $X^{\mathcal
{V}\setminus
A}(t)$ have common jump times contradicting that $X(t)$ is composable, or
they constitute 2 separate SKMs which then require a new method for their
individual analysis.

The kinetic independence graph of an SKM is defined as follows.
\begin{definition}
\label{LIGnew}The directed graph $G$ with vertex set $\mathcal{V}$ is the
\textit{kinetic independence graph} (\textit{KIG}) of the SKM $[N,S,\mathsf{P}]$ if
and only if
%
%
\begin{equation}\label{Newdef}
\mathrm{pa}(k)=R[\Delta(k)]\setminus\{k\}\qquad \forall k\in\mathcal{V},
\end{equation}
where $\mathrm{pa}(k)=\{i\in\mathcal{V}|i\rightarrow k\}$ is the set of parents of
vertex $k$, and $R[\Delta(k)]:=\bigcup_{m \in\Delta(k)}R[m]$ is the
set of reactants of all
reactions that change species $k$.
\end{definition}

Since only partial information about the SKM is required for
construction of the KIG, the necessary information is currently
available for many biochemical
reaction networks. For each $m\in\mathcal{M}$, it is required to know
the reactants $R[m]$, and the species (reactants and products) changed by
the reaction, that is, $\{i\in\mathcal{V}|S_{im}\neq0\}$. Full
knowledge of the
stoichiometric matrix $S$ is neither necessary nor sufficient for
construction of the KIG. Note that the possible presence of a catalyst
among the reactants
$R[m]$ implies that $R[m]$ cannot be reliably reconstructed from $S$.
No knowledge of the rate
parameters $c_{m}$ is required for construction of the KIG, which is
important since their measurement is difficult experimentally.

A comment will be useful at this juncture on the treatment of
measurability considerations in the paper. While the treatment is fully
rigorous, it is appreciated that some readers will be more concerned
with application of the paper's results. Proofs requiring a
measure-theoretic approach have therefore been placed in Appendices~\ref{appA}
and \ref{Prooffs}. Note that a statement such as the one that $\lambda_{m}(t)$ is
(as it must be) measurable $\mathcal{F}_{t}^{R[m]}$ implies that the
realized value of the r.v. $\lambda_{m}(t)$ may be ``computed'' from the
sample path of the subprocess for $R[m]$ over the interval $[0,t]$.

The motivation for Definition \ref{LIGnew} of the KIG of an SKM is that
the local evolution of species $k$ depends only on the stochastic rate
of reactions that change the number of molecules (the level) of $k$,
which in turn depend only on the levels of their reactants. To make
this exact, the concept of local independence \cite{Didelez2008} is
needed. Let $A,B\subset\mathcal{V}$. We will say that $N^{B}$ is
\textit{locally independent} of $N^{A}$
(given $N^{\mathcal{V}\setminus A}$) if and only if the $\mathcal
{F}_{t}^{N}$-intensity, $\lambda^{B}(t)$, is measurable $\mathcal
{F}_{t}^{\mathcal{V}\setminus A}$ for all $t$---that is, the internal
history of $X_{t}^{A}$ is irrelevant for the $\mathcal
{F}_{t}^{N}$-intensity of the species in $B$. Only intensities of
subprocesses conditional on the history of the whole system, $\mathcal
{F}_{t}^{N}$, are considered here (as opposed to $\mathcal
{G}_{t}$-intensities where $\mathcal{G}_{t}\subset\mathcal{F}_{t}^{N}$).

As a consequence of Definition \ref{LIGnew}, one can read off from the
KIG, for
any collection of vertices $B$, those subprocesses with respect
to which $N^{B}$ is locally independent, that is, which are irrelevant for
the instantaneous evolution of $B$. Denote the closure of $B$ by
$\mathrm{cl}(B):=\mathrm{pa}(B)\cup B$.
\begin{proposition}
\label{local}Let $G$ be the KIG of the SKM $[N,S,\mathsf{P}]$ and let
$A,B\subset\mathcal{V}$. Then the $%
\mathcal{F}_{t}^{N}$-intensity $\lambda^{B}(t)$~is measurable $\mathcal
{F}%
_{t}^{\mathrm{cl}(B)}$ for all $t$, that is, $N^{B}$ is \textit{locally independent}
of $N^{\mathcal{V}\setminus \mathrm{cl}(B)}$ (given $N^{\mathrm{cl}(B)}$). Suppose that
$A\cap \mathrm{cl}(B)=\varnothing$. Then $\lambda^{B}(t)$ is measurable $\mathcal
{F}%
_{t}^{V\setminus A}$.
\end{proposition}
\begin{pf}
By (\ref{SupInt}), each intensity $\lambda_{e}^{B}(t)$ is
measurable $%
\mathcal{F}_{t}^{R[\Delta(B)]}$ because $\lambda
_{m}(t)=c_{m}g_{m}\{X^{R[m]}(t-)\}$ is measurable $\mathcal{F}%
_{t}^{X^{R[\Delta(B)]}}=\mathcal{F}_{t}^{R[\Delta(B)]}$ $\forall m\in
\Delta(B)$, recalling that $R[\Delta(B)]=\bigcup_{m\in\Delta
(B)}R[m]$%
. Since $\mathrm{pa}(B)=\{\bigcup_{k\in B}\mathrm{pa}(k)\}\setminus\{B\}=\{\bigcup
_{k\in
B}R[\Delta(k)]\}\setminus\{B\}=\{R[\Delta(B)]\}\setminus\{B\}$, it
follows that $R[\Delta(B)]\subseteq \mathrm{cl}(B)$, and hence $\mathcal{F}%
_{t}^{R[\Delta(B)]}\subseteq\mathcal{F}_{t}^{\mathrm{cl}(B)}$ by Lemma \ref{Filts}.
Thus, each intensity $\lambda_{e}^{B}(t)$ is measurable $\mathcal{F}%
_{t}^{\mathrm{cl}(B)}$, and the remainder of the proposition follows immediately.
\end{pf}

Proposition \ref{local} accords with chemical intuition.
Given the internal history
of $X^{\mathrm{cl}(B)}$ at time $t$, the levels of the species $R[\Delta(B)]$
just prior to $t$ are ``known.'' These are exactly the species levels
that determine the local dynamics of $B$ since, as reactants, they
determine the rate of all reactions that change the concentrations of
$B$. Therefore, any further information about species histories,
including the internal history of $N^{\mathcal{V}\setminus \mathrm{cl}(B)}$, is
irrelevant for the local dynamics of $B$.

Notice that loops, that is edges of the type $k\rightarrow k$ are by
definition not included in the KIG, even though $k$ may well be in
$R[\Delta(k)]$. For this reason, one cannot assert in
Proposition \ref{local} that $\lambda^{B}(t)$ is
measurable $\mathcal{F}_{t}^{\mathrm{pa}(B)}$, but rather that it is measurable
$\mathcal{F}_{t}^{\mathrm{cl}(B)}$. More generally, a particular SKM may imply further
local independences of $\lambda^{B}(t)$ than those encoded by the
KIG---for example, due to a
deterministic relationship between two subsets of species arising from a
chemical conservation relation---but this level of knowledge
about the SKM is not assumed in constructing the KIG.

Graphical separations in the undirected version of the KIG, written $G^{
\thicksim}$, are central in what follows. Diagrammatically, $G^{%
\thicksim}$ is the undirected graph obtained from $G$ by substituting
lines for arrows.
Let $A,B,D\subset\mathcal{V}$.
The notation $A\perp_{G^{\thicksim}}B|D$ stands for the \textit{graphical
separation} of $A$ from $B$ by $D$, that is, the property that every sequence
of edges (or path) in $G^{\thicksim}$ that begins with some $a\in A$
and (without any repetition of vertices) ends with some $b\in B$,
includes a
vertex in~$D$. With $[A,B,D]$ a partition of $\mathcal{V}$, such a
separation in $G^{\thicksim}$ is equivalent to the nonexistence of
$(a\in A,b\in B)$ such that there is an edge $a\rightarrow b$ or an
edge $b\rightarrow a$ in $G$. This graphical separation implies the
following mutual local independence property.
\begin{proposition}
\label{bdli}Let $G$ be the KIG of an SKM $[N,S,\mathsf{P}]$, and let $[A,B,D]
$ be a partition of $\mathcal{V}$. If $A\bot_{G^{\thicksim}}B|D$, then
$%
N^{B}$ is \textit{locally independent} of $N^{A}$ (given $N^{B \cup
D}$) and $N^{A}$ is \textit{%
locally independent} of $N^{B}$ (given $N^{A \cup D}$) or,
equivalently, $\lambda^{B}(t)$ is
measurable $\mathcal{F}_{t}^{B\cup D}$ and $\lambda^{A}(t)$ is
measurable $%
\mathcal{F}_{t}^{A\cup D}$.
\end{proposition}
\begin{pf}
$A\bot_{G^{\thicksim}}B|D$, if and only if $B\cap \mathrm{cl}(A)=A\cap
\mathrm{cl}(B)=\varnothing$.
The result then follows directly from Proposition \ref{local}.
\end{pf}

Note that it follows from the definition of the KIG that the graphical
separation in Proposition \ref{bdli} is equivalent to the chemical property
$A\cap R[\Delta(B)]=B\cap R[\Delta(A)]=\varnothing$. That is, $A$ does
not participate as a reactant in any reaction that changes $B$,
and vice versa. Therefore, for example, $R[\Delta(B)]
\subseteq B \cup D$---hence, given the levels of $B$ and $D$ (which
fully determine the rate of reactions that change~$B$), the levels of
$A$ are irrelevant for the instantaneous evolution of $B$, and
vice versa. Section \ref{sec4} will establish that under weak
regularity conditions on the SKM, the
separation $A\bot_{G^{\thicksim}}B|D $ in $G^{\thicksim}$ implies not
only mutual local independence but also \textit{global} conditional
independence of the internal histories of $A$ and $B$ given a history
of $D$.

\begin{figure}

\includegraphics{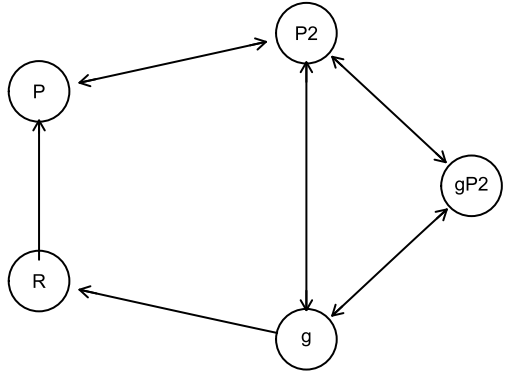}

\caption{Kinetic independence graph of the SKM in
Example \protect\ref{Example1_WellBeh}.}\label{FigExample1_WellBeh}
\end{figure}

As an illustration of the concepts discussed so far, consider again the
SKM of Example \ref{Example1_WellBeh}. The corresponding KIG is shown
in Figure \ref{FigExample1_WellBeh}.
Note the presence of cycles in the KIG, including $g\rightarrow
R\rightarrow
P\rightarrow P_{2}\rightarrow g$ which might be termed the ``negative
feedback cycle.'' Clearly, $\{P,R\}\bot_{G^{\thicksim}}$ $\{gP_{2}\}|%
\{g,P_{2}\}$. Let $D:=\{g,P_{2}\}$. Notice that, according to
Definition \ref{NADeltaA},
SKM subprocesses are given by $N^{gP_{2}}=[N_{b},N_{ub}]^{\prime}$ and
$%
N^{D}=[N_{d},N_{rd},N_{b},N_{ub}]^{\prime}$. Hence, $\mathcal{F}%
_{t}^{gP_{2}}\subset\mathcal{F}_{t}^{D}$, and the global independence
$\mathcal{F}_{t}^{P,R}\indep{\mathcal{F}%
_{t}^{gP_{2}}}|\mathcal{F}_{t}^{D}$ holds immediately in this case.

%

Anticipating the problem of how to modularize SKMs to be tackled in
Section~\ref{sec5}, a modularization suggested for the SKM of
Example \ref{Example1_WellBeh} by its dynamic independence properties is the
modularization $[\{P,R,g,P_{2}\}$, $\{gP_{2},g,P_{2}\}]$. The 2 module
``residuals'' are given by $\{P,R\}$ and $\{gP_{2}\}$. Each module
residual is locally independent of the other given that module's
internal history. Furthermore, the 2 modules are conditionally
independent given the history of their intersection, $\mathcal
{F}_{t}^{g,P_{2}}$. In fact, these 2 modules correspond to the maximal
prime subgraphs of $G^{\thicksim}$ for this example (see Definition
\ref{MPD}). The graphical
methods for identifying SKM modularizations in Section \ref{sec5}
are, broadly speaking,
also based around the maximal prime decomposition of the undirected KIG.

\section{Global dynamic independence}
\label{sec4}

This section will present the theorems establishing that for
a partition $[A,B,D]$ of
$\mathcal{V}$, the separation of $A$ from $B$ by $D$ in the undirected
version of the KIG implies the global dynamic independence $\mathcal{F}%
_{t}^{A}\indep\mathcal{F}_{t}^{B}|%
\mathcal{F}_{t}^{D^{\ast}}$ for all $t\geq0$, under the probability
measure of the SKM, $\mathsf{P}$.

The history of $X^{D}$ given by $\mathcal{F%
}_{t}^{D^{\ast}}$ is defined formally below. Heuristically, $\mathcal
{F}_{t}^{D^{\ast}}$~includes at time $t$ the internal history of the
jump process for the species in $D$ $(\mathcal{F}%
_{t}^{D} \subseteq\mathcal{F}_{t}^{D^{\ast}})$, and also for every
jump time of $X^{D}$ always contains the
``information'' whether some species in $A$ jumped, or some species in $B$
jumped, or some species in both $A$ and $B$ jumped, but not (necessarily)
the particular species involved. The main proofs of the theorems are
quite involved and are given in Appendix \ref{appA}. Readers less concerned with
technical details will find an outline of the argument and intuitions
for important aspects of the proofs in this section of the paper. It is
worth defining here explicitly what is meant by the conditional
independence of $\sigma$-fields \cite{DMeyer78,FMR90}.
\begin{definition}
\label{CI}Let $(\Omega,\mathcal{F,}\mathsf{P})$ be\vspace*{1pt} an arbitrary probability
space and suppose we have 3 sub-$\sigma$-fields $\mathcal
{F}^{1},\mathcal{F}%
^{2},\mathcal{F}^{3}\subseteq\mathcal{F}$. We say that $\mathcal{F}^{1}$
and $\mathcal{F}^{2}$are independent conditionally on $\mathcal{F}^{3}$ and
write $\mathcal{F}^{1}\indep\mathcal{F}%
^{2}|\mathcal{F}^{3};\mathsf{P}$ if and only if
\[
\mathsf{E}[Z_{1}|\mathcal{F}^{2}\vee\mathcal{F}^{3}]=\mathsf{E}[Z_{1}|%
\mathcal{F}^{3}]
\]
for all nonnegative random variables $Z_{1}$ that are measurable
$\mathcal{F%
}^{1}$. The notation $\mathcal{F}^{2}\vee\mathcal{F}^{3}$ stands for
the smallest $\sigma$-field
containing both $\mathcal{F}^{2}$ and $\mathcal{F}^{3}$. The
relationship is symmetric, that is,
$\mathcal{F}^{1}\indep\mathcal{F}^{2}|\mathcal{F}^{3};\mathsf{%
P}\Leftrightarrow\mathcal{F}^{2}\indep{%
\mathcal{F}^{1}}|\mathcal{F}^{3};\mathsf{P}$.
\end{definition}

Thus, the global dynamic independence statement $\mathcal{F}%
_{t}^{A}\indep{\mathcal{F}_{t}^{B}}|%
\mathcal{F}_{t}^{D^{\ast}}$ can be understood as follows: the
expectation\vspace*{1pt} of (suitably measurable) mappings from sample paths of
$N^{A}$ (resp., $N^{B}$) over $(0,t]$ to $\mathbb{R}$, conditional on
the history $\mathcal{F}_{t}^{D^{\ast}}$, are
unchanged when the conditioning $\sigma$-field also includes the
internal history of $N^{B}$ (resp., $N^{A}$). Roughly speaking, and over
any time interval $(0,t]$, all ``information'' about the dynamic
evolution of $B$ is irrelevant for the dynamic evolution of $A$, given
the ``information'' in $\mathcal{F}_{t}^{D^{\ast}}$ (and vice versa).

First, an outline of the logic of the argument of this section is presented,
before going on to state the main theorems.

\subsection{Preliminaries and outline of argument}\label{sec41}

The following lemma is central to the method used. Although closely related
to a result in \cite{FMR90}, I am not aware of its statement and proof
elsewhere.
\begin{lemma}
\label{CIDom}Let $\mathsf{P},\mathsf{\tilde{P}}$ be probability
measures on
an arbitrary measurable space, $(\Omega,\mathcal{F})$, such that
$\mathsf{P}%
\ll\mathsf{\tilde{P}}$. Consider any 3 sub-$\sigma$-fields $\mathcal
{F}%
^{1},\mathcal{F}^{2},\mathcal{F}^{3}\subseteq\mathcal{F}$ satisfying the
conditional\vspace*{1pt} independence $\mathcal{F}^{1}\indep
{\mathcal{F}^{2}}|\mathcal{F}^{3};\mathsf{\tilde{P}}$ under the
dominating measure $\mathsf{\tilde{P}}$. Denote by $\mathcal{L}_{123}$ a
Radon--Nikodym derivative, $(\mathsf{dP}_{t}/\mathsf{d\tilde{P}}_{t})|_{
\mathcal{F}^{1}\vee\mathcal{F}^{2}\vee\mathcal{F}^{3}}$.

Then the following condition implies that the conditional independence $
\mathcal{F}^{1}\indep{\mathcal{F}%
^{2}}|\mathcal{F}^{3}$ holds also under $\mathsf{P}$:
\[
\mathcal{L}_{123}=\psi_{13}\psi_{23},
\]
where $\psi_{i3}$ is a nonnegative, $\mathcal{F}^{i}\vee\mathcal
{F}^{3}$%
-measurable random variable for $i\in\{1,2\}$.
\end{lemma}

Proof of Lemma \ref{CIDom} is given in Appendix \ref{Prooffs}.

Since $\mathcal{F}_{t}^{D^{\ast}}$ is a history of $N^{D}$ (i.e.,
$\mathcal{F}_{t}^{D} \subseteq\mathcal{F}_{t}^{D^{\ast}}$), it
follows from Lemma \ref{Filts} that $\mathcal{F}_{t}^{A}\vee\mathcal
{F}_{t}^{B}\vee\mathcal{F}_{t}^{D^{\ast}}=\mathcal{F}_{t}^{X}=\mathcal
{F}_{t}^{N}$. A likelihood process $\mathcal{L}_{t}:=(\mathsf{dP}_{t}/%
\mathsf{d\tilde{P}}_{t})|_{\mathcal{F}_{t}^{N}}$ is thus required in order
to apply Lemma \ref{CIDom} to the 3 $\sigma$-fields $\mathcal
{F}_{t}^{A},%
\mathcal{F}_{t}^{B},\mathcal{F}_{t}^{D^{\ast}}$. Given its importance here,
we restate for an SKM the following likelihood result from the counting
process literature (see, e.g., \cite{Jacobsen82},
Theorem 4.1, page 74). The proof is
omitted since it is well known.
\begin{lemma}
\label{Likl}Let $[N,S,\mathsf{P}]$ be an SKM as in Definition
\ref{SKM}, and
let $[N,\mathsf{\tilde{P}}]$ be the $M$-variate Poisson process with
intensities $\mathbf{1}_{M}=(1,\ldots,1)^{\prime}$. Then, for every
$t\geq0$,
$\mathsf{P}_{t}\ll\mathsf{\tilde{P}}_{t}$ and a Radon--Nikodym
derivative is
given by
%
%
\begin{equation} \label{RNderiv}
\frac{\mathsf{dP}_{t}}{\mathsf{d\tilde{P}}_{t}}\bigg|_{\mathcal{F}%
_{t}^{N}}=\prod_{m=1}^{M} \biggl\{ \prod_{T_{s}^{m}\leq t}\lambda
_{m}(X_{T_{s}^{m}-}) \biggr\} \exp\biggl\{ t-\int_{0}^{t}\lambda
_{m}(X_{u-})\,du \biggr\}.
\end{equation}
Note that the counting processes $[N_{m};m=1,\ldots,M]$ are independent
under $%
\mathsf{\tilde{P}}$ (see, e.g., \cite{Jacobsen2006}, Proposition
4.7.2), and
hence the $\sigma$-fields $[\mathcal{F}_{t}^{m};m=1,\ldots,M]$ are independent
under $\mathsf{\tilde{P}}$ and $\mathsf{\tilde{P}}_{t}$ for all $t\geq0$.
\end{lemma}

Of considerable importance here will be the fact that, under the
dominating measure $\mathsf{\tilde{P}}$, the counting processes $%
[N_{m};m=1,\ldots,M]$ are independent. Of course, two or more of the
subprocesses $[N^{A},N^{B},N^{D}]$ may have jump times in common as the
result of reactions that simultaneously change several of the species
sets $[A,B,D]$. However, denoting by $\Delta D_{D}$ the reactions that
change $D$ alone, the reaction set $\mathcal{M}$ can be partitioned as
$[\Delta(A),\Delta(B)\setminus\Delta(A),\Delta D_{D}]$.
The independence of the reaction counting processes then implies that
$\mathcal{F}%
_{t}^{\Delta(A)}\indep{\mathcal{F}%
_{t}^{\Delta(B)\setminus\Delta(A)}}|\mathcal{F}_{t}^{\Delta D_{D}};%
\mathsf{\tilde{P}}_{t}$, which is a point of departure for proving
Theorem \ref{Main1} below.

To apply Lemma \ref{CIDom}, we first establish that if $A\bot
_{G^{\thicksim}}B|D$, then $\mathcal{F}_{t}^{A}\indep{%
\mathcal{F}_{t}^{B}}|\mathcal{F}_{t}^{D^{\ast}}$ under the dominating
measure $\mathsf{\tilde{P}}$ (see Theorem \ref{Main1}). We then show
that the factorisation $\mathcal{L}_{t}=\psi
_{AD^{\ast},t}\psi_{BD^{\ast},t}$ holds with, for example, $\psi
_{AD^{\ast},t}$
an $\mathcal{F}_{t}^{A}\vee\mathcal{F}_{t}^{D^{\ast}}$-measurable r.v.
(see Theorem \ref{Main2}). We are thus able to conclude that if $A\bot
_{G^{%
\thicksim}}B|D$, then the SKM must satisfy $\mathcal{F}_{t}^{A}\indep
{\mathcal{F}_{t}^{B}}|\mathcal{F}%
_{t}^{D\ast};\mathsf{P}$ (see Corollary \ref{Main3}).

Definition of the filtration $\{\mathcal{F}_{t}^{D^{\ast}}\}$ is needed.
This is best understood as the internal history of a particular MVCP, $%
N^{D^{\ast}}(t)$, defined now below.
\begin{definition}
\label{NDStar}Let $[A,B,D]$ be a partition of $\mathcal{V}$, the
species set of an SKM. Define $\Delta D_{A}:=\{\Delta(D)\cap\Delta
(A)\}\setminus\Delta(B)$, the set of reactions that change $D$ and $A$,
but not $B$. Similarly, define $\Delta D_{AB}:=\Delta(D)\cap\Delta
(A)\cap
\Delta(B)$, $\Delta D_{B}:=\{\Delta(D)\cap\Delta(B)\}\setminus
\Delta
(A)$, and those that change $D$ alone by $\Delta D_{D}:=\Delta
(D)\setminus
\{\Delta(A)\cup\Delta(B)\}$. Then $[\Delta D_{A},\Delta D_{AB},\Delta
D_{B},\Delta D_{D}]$ is a partition of $\Delta(D)$, the reactions that
change~$D$.\vspace*{1pt}

The MVCP $N^{D^{\ast}}(t)$ is constructed by taking
each element, $\Delta D_{\bullet}$, of this partition in turn and summing
over counting processes for reactions in $\Delta D_{\bullet}$ alone that
result in identical changes to $D$---that is, by applying to each
$\Delta
D_{\bullet}$ the equivalence relation
\[
m\sim m^{\prime}\quad\Leftrightarrow\quad S_{m}^{D}=S_{m^{\prime}}^{D},\qquad
m,m^{\prime} \in\Delta D_{\bullet}.
\]
The resultant MVCP, denoted $N_{\bullet}^{D}(t)$, is given by
\[
N_{\bullet}^{D}(t):= \biggl\{ \sum_{m\in\mathcal{M}_{e}(\Delta
D_{\bullet
})}N_{m}(t) \biggr\}_{e=1,\ldots,|\mathcal{M}(\Delta D_{\bullet})|},
\]
which differs from the SKM subprocess $N^{D}(t)$ of Definition
\ref{NADeltaA} in that the partition $\mathcal{M}(\Delta D_{\bullet})$ is used
in place of $\mathcal{M}(\Delta(D))$. Then define
\[
N^{D^{\ast}}(t):=[N_{A}^{D}(t),N_{AB}^{D}(t),N_{B}^{D}(t),N_{D}^{D}(t)].
\]
Here $\Delta D_{\bullet}$ is empty, the relevant component of
$N^{D^{\ast
}}(t)$ is set equal to zero $\forall t\geq0$. The corresponding internal
histories of $N_{{\bullet}}^{D}(t)$ are written $\{\mathcal{F}%
_{A}^{D}(t),\break\mathcal{F}_{AB}^{D}(t),\mathcal{F}_{B}^{D}(t)$, $\mathcal{F}%
_{D}^{D}(t)\}$.
\end{definition}

Thus, as stated previously, at time $t$ $\mathcal
{F}_{t}^{D^{\ast}}$ includes the
internal history of the jump process for the species in $D$, $\mathcal
{F}%
_{t}^{D}$, and also for every jump time of $X^{D}$ always contains the
``information'' whether some species in $A$ jumped, or some species in $B$
jumped, or some species in both $A$ and $B$ jumped, but not (necessarily)
the particular species involved. An alternative formulation of $\mathcal
{F}_{t}^{D^{\ast}}$
would be as the internal history at $t$ of the marked point process $\{
T_{s}^{D},\tilde{Z}_{s}^{D}\}$, where the marks $\tilde{Z}_{s}^{D}$
give not only the value of the jump in $D$ but also an indicator of
which element, $\Delta D_{\bullet}$, the reaction causing the jump in
$D$ belongs to.

In some applications, it may be more
convenient or practical to use only internal histories. Section \ref{sec43}
will thus provide a rigorous and intuitive means of comparing
$N^{D^{\ast
}}(t)$ and $N^{D}(t)$---through comparison of the corresponding partitions
of the reactions in $\Delta(D)$---and state a property that is easily
checked for a given SKM under which $\mathcal{F}_{t}^{D^{\ast
}}=\mathcal{F}%
_{t}^{D}$.

The following conditions on the SKM are used in the statement of the results
of this section. Both Theorems \ref{Main1} and \ref{Main2} assume that the
SKM is standard, which imposes the following very weak regularity
conditions.
\begin{definition}
\label{stdSKM}An SKM $[N,S,\mathsf{P}]$ is a \textit{standard SKM} if it
satisfies all of the following: (i) every reaction changes at least 1
species, that is, $S_{m}\neq\mathbf{0}$ $\forall m\in\{1,\ldots,M\}$;
(ii) every species in $\mathcal{V}$ is changed by at least one reaction,
that is, the row $S_{k\bullet}\neq\mathbf{0}$ $\forall k\in\mathcal
{V};$ (iii) if a
zeroth order reaction $\tilde{m}$ is included (i.e., $R[%
\tilde{m}]=\varnothing$) then it has only 1 product; (iv) for all $m$,
if $%
|R[m]|=1$ then $|R^{\ast}[m]|=1$ and if $|R[m]|>1$ then $%
|\{R[m]\}\setminus\{R^{\ast}[m]\}|\leq1$, where $R^{\ast}[m]=\{i\in
R[m]|S_{im}\neq0\}$ are the reactants changed by $m$.
\end{definition}

The first condition of Definition \ref{stdSKM} is obvious. The
second does not preclude an effect of the concentration of species that are
constant over time (via the functions $g_{m}$ or the constants
$c_{m}$). The
third is just a convention. The fourth ensures that if $R[m]\neq
\varnothing$ then the reaction has at least 1 reactant that is changed
and at
most 1 reactant that is not changed. It allows for the inclusion of a
reaction with an unchanged reactant where this simplifies the SKM, for example,
where that reactant acts as a catalyst [although the reaction could be
broken down into several reactions not requiring condition (iv) if
desired].

Theorems \ref{Main1} and \ref{Main2} also require that the following
condition holds for $\Gamma=\Delta(A)\cap\Delta(B)$. If $\Delta
(A)\cap
\Delta(B)=\varnothing$, as sometimes happens, then the condition is trivial
and always satisfied.
\begin{condition}
\label{IdentGamma} Let $[N,S,\mathsf{P}]$ be a standard
SKM. A subset of reactions $\Gamma$, $\varnothing\subseteq
\Gamma
\subseteq\mathcal{M}$, is said to be identified by consumption of reactants
if and only if:

(i) For all $m\in\Gamma$, $S_{im}\leq0$ $\forall i\in R[m]$ (hence, $
S_{km}<0$ for some $k\in R[m]$ provided that $R[m]\neq\varnothing$),
and $%
S_{im}\geq0$ $\forall i\in P[m]$; and (ii) $\nexists m,\tilde{m}\in
\Gamma$
$(m\neq\tilde{m})$ such that $S_{m}^{-}=S_{\tilde{m}}^{-}$, where $%
S_{m}^{-} $ denotes the vector formed by setting all positive elements
of $%
S_{m}$ to zero.
\end{condition}
\begin{remark}
\label{OnIdentGam}
Condition \ref{IdentGamma} implies that no 2 reactions in $%
\Gamma$ change reactants identically, hence the reactions in $\Gamma$ are
identified uniquely by their consumption of reactants. Condition
\ref{IdentGamma} will be satisfied with $\Gamma=\mathcal{M}$ by most SKMs of
interest, possibly after explicit inclusion of enzymes in reaction
mechanisms. Although autocatalytic reactions such as
$X_{j}+X_{k}\rightarrow
^{m}2X_{k}$ and its reverse violate condition~(i), these could be
accommodated by instead including a more detailed mechanism, for
example, $%
X_{j}+X_{k}\rightarrow^{m_{1}}X_{j}X_{k}$ and $X_{j}X_{k}\rightarrow
^{m_{2}}2X_{k}$.

An alternative approach would be to work with $G_{f}^{\thicksim}$, the
graph obtained from the undirected version of the KIG by adding an edge
$%
j\sim k$ whenever $[S_{jm}>0$ and $S_{km}>0$ for some $m\in\mathcal{M}]$
and $j\nsim k$ in $G^{\thicksim}$. The graph $G_{f}^{\thicksim}$ might
be termed the fraternized (as distinct from moralized) version of the
KIG. The separation $A\bot_{G_{f}^{%
\thicksim}}B|D$ implies that $\Delta(A)\cap\Delta(B)=\varnothing$
[since $A\bot_{G^{\thicksim}}B|D$ and hence $R[m]\subseteq D$ for any $
m\in\Delta(A)\cap\Delta(B)$]. Therefore, if $A\bot_{G^{\thicksim%
}}B|D $ is replaced by $A\bot_{G_{f}^{\thicksim}}B|D$, Condition
\ref{IdentGamma} for $\Gamma=\Delta(A)\cap\Delta(B)$ can be dropped
from the
statements of Theorems \ref{Main1} and \ref{Main2}, and from that of
Corollary \ref{Main3}.
\end{remark}

\subsection{Global independence theorems}\label{sec42}

We are now in a position to state the main results of Section \ref{sec4}
of the paper. Theorem \ref{Main1} is concerned with global dynamic independence
under $\mathsf{\tilde{P}}$, the law of the $M$-variate Poisson process
(see Lemma \ref{Likl}).
\begin{theorem}
\label{Main1}Let $G$ be the KIG of a standard SKM, $[N,S,\mathsf{P}]$, and
let $[A,B,D]$ be a partition of $\mathcal{V}$. Suppose also that Condition
\ref{IdentGamma} holds for $\Gamma=\Delta(A)\cap\Delta(B)$ (where $%
\Gamma$ is possibly empty, in which case the condition is trivial).
Then $%
A\bot_{G^{\thicksim}}B|D$ implies that $\mathcal{F}_{t}^{A}\indep
\mathcal{F}_{t}^{B}|\mathcal{F}%
_{t}^{D^{\ast}};\mathsf{\tilde{P}}_{t}$, where $\{\mathcal
{F}_{t}^{D^{\ast
}}\}$ is the natural filtration of $N^{D^{\ast}}(t)$, $N^{D^{\ast
}}(t)$ is
given by Definition \ref{NDStar}, and $\mathsf{\tilde{P}}$ is the law
of the
$M$-variate Poisson process in Lemma \ref{Likl}.
\end{theorem}

We provide here a somewhat heuristic discussion of this result, a
rigorous treatment
being given in Appendix \ref{appA1}. The argument can be broken down into four steps.

First, the reaction counting processes $[N_{m};m=1,\ldots,M]$ are independent
under $\mathsf{\tilde{P}}$. Therefore,
%
%
\begin{equation} \label{one}
\mathcal{F}_{t}^{\Delta(A)}\indep{%
\mathcal{F}_{t}^{\Delta(B)\setminus\Delta(A)}}|\mathcal
{F}_{t}^{\Delta
D_{D}};\mathsf{\tilde{P}}_{t},
\end{equation}
since $[\Delta(A),\Delta(B)\setminus\Delta(A),\Delta D_{D}]$ is a
partition of the reaction set $\mathcal{M}$. Equation (\ref{one})
holds because the three MVCPs associated with each element of the
partition are
(unconditionally) independent.

Second, consider again Definition \ref{NDStar} for $N^{D^{\ast}}(t)=[\{
N_{A}^{D}(t),N_{AB}^{D}(t)\}$,\break $\{N_{B}^{D}(t)\},\{N_{D}^{D}(t)\}]$.
The internal history of the first component MVCP in curly parentheses
must be contained in the internal
history of $N^{\Delta(A)}(t)$. [All the reactions involved
in that component change $A$ and hence the sample path of $N^{\Delta
(A)}(t)$ implies that of the first component.] Similarly, the internal
history of the second component of $%
N^{D^{\ast}}(t)$ in curly parentheses must be contained in that of
$N^{\Delta(B)\setminus
\Delta(A)}(t)$. The internal history of the third component is equal
to $\mathcal{F}%
_{t}^{\Delta D_{D}}$. Combining the internal histories of
these 3 components making up $N^{D^{\ast}}(t)$ must give $\mathcal{F}%
_{t}^{D^{\ast}}$. Therefore, the internal histories of the first 2
components can be used to expand the conditioning information in
(\ref{one}) to give
%
%
\begin{equation} \label{two}
\mathcal{F}_{t}^{\Delta(A)}\indep{%
\mathcal{F}_{t}^{\Delta(B)\setminus\Delta(A)}}|\mathcal
{F}_{t}^{D^{\ast
}};\mathsf{\tilde{P}}_{t}.
\end{equation}

Third, establishing the property $\mathcal{F}_{t}^{\Delta(A)}\indep
{\mathcal{F}_{t}^{\Delta(B)}}|\mathcal{F}%
_{t}^{D^{\ast}};\mathsf{\tilde{P}}_{t}$ implies the global dynamic
independence in Theorem \ref{Main1} since the internal history of the subprocess
for $A$ must be contained in that of $N^{\Delta(A)}$ [since the sample path
of $N^{\Delta(A)}(t)$ obviously implies that of $N^{A}(t)$], and similarly
for $B$. This property in turn follows by showing that the internal history
of $N^{\Delta(A)\cap\Delta(B)}(t)$ is contained in $\mathcal{F}%
_{t}^{D^{\ast}}$. The second $\sigma$-field in (\ref{two}) can
then be expanded to include $\mathcal{F}_{t}^{\Delta(A)\cap\Delta(B)}
$. Combining the internal histories $\mathcal{F}_{t}^{\Delta(A)\cap
\Delta
(B)}$ and ${\mathcal{F}_{t}^{\Delta(B)\setminus\Delta(A)}}$ in this way
gives ${\mathcal{F}_{t}^{\Delta(B)}}$.

Finally, $\mathcal{F}_{t}^{\Delta(A)\cap\Delta(B)}\subseteq\mathcal
{F}%
_{t}^{D^{\ast}}$is a direct consequence of the fact that $\Delta
D_{AB}=\Delta(A)\cap\Delta(B)$---that is, all reactions that change
$A$ and
$B$ also change $D$---and that the reactions in $\Delta D_{AB}$ change $D$
uniquely (among themselves). These properties of $\Delta D_{AB}$ depend
crucially on the graphical separation $A\bot_{G^{\thicksim}}B|D$ and
also on Condition \ref{IdentGamma} holding for $\Delta D_{AB}$ (under the
conditions of Theorem \ref{Main1}). The separation ensures that for any $m\in
\Delta(A)\cap\Delta(B)$, the reactants of $m$ are all in $D$
(otherwise, we would have
either $A\rightarrow B$ or $B\rightarrow A$ in the KIG) and hence also
$m \in\Delta(D)$. Condition \ref{IdentGamma} ensures that the members
of $\Delta D_{AB}$ are identified by
consumption of reactants, hence the reactions in $\Delta D_{AB}$ must change
$D$ uniquely (among themselves) and $\mathcal{F}_{t}^{\Delta
(A)\cap\Delta(B)}=\mathcal{F}_{AB}^{D}(t)$. Therefore,
$\mathcal{F}_{t}^{\Delta(A)\cap\Delta(B)}\subseteq\mathcal{F}%
_{t}^{D^{\ast}}$, since $N_{AB}^{D}(t)$ is a component of $N^{D^{\ast}}(t)$.

We now turn to consider global dynamic independence
under $\mathsf{P}$, the law of the SKM.
\begin{theorem}
\label{Main2}Let $G$ be the KIG of a standard SKM, $[N,S,\mathsf{P}]$, and
let $[A,B,D]$ be a partition of $\mathcal{V}$. Suppose also that Condition
\ref{IdentGamma} holds for $\Gamma=\Delta(A)\cap\Delta(B)$ (where $%
\Gamma$ is possibly empty, in which case the condition is trivial).
Then $%
A\bot_{G^{\thicksim}}B|D$ implies that
\[
\mathcal{L}_{t}:=(\mathsf{dP}_{t}/\mathsf{d\tilde{P}}_{t})|_{\mathcal{F}
_{t}^{N}}=\psi_{AD^{\ast},t}\cdot\psi_{BD^{\ast},t},\qquad t\geq0,
\]
where $\psi_{iD^{\ast},t}$ is a nonnegative, $\mathcal{F}_{t}^{i}\vee
\mathcal{F}_{t}^{D^{\ast}}$-measurable random variable for $i\in\{
A,B\}$,
and $\{\mathcal{F}_{t}^{D^{\ast}}\}$ is the natural filtration of $%
N^{D^{\ast}}(t)$.
\end{theorem}

Taking the logarithm of the likelihood in (\ref{RNderiv}) yields
%
%
\begin{eqnarray} \label{three}
\log\mathcal{L}_{t}&=&\sum_{m\in\mathcal{M}}\biggl[t-\int_{0}^{t}\lambda
_{m}(u)\,du+\sum_{s\geq1}1(T_{s}^{m}\leq t)\log(\lambda
_{m}(T_{s}^{m}))\biggr]\nonumber\\[-8pt]\\[-8pt]
:\!&=&\sum_{m\in\mathcal{M}}l_{m}(t).\nonumber
\end{eqnarray}
Theorem \ref{Main2} may be established by showing that, $\forall m\in
\mathcal{M}$, $l_{m}(t)$ is measurable either $\mathcal
{F}_{t}^{AD^{\ast}}:=%
\mathcal{F}_{t}^{A}\vee\mathcal{F}_{t}^{D^{\ast}}$ or $\mathcal{F}%
_{t}^{BD^{\ast}}:=\mathcal{F}_{t}^{B}\vee\mathcal{F}_{t}^{D^{\ast
}}$. We
explain here how $l_{m}(t)$ may be
computed ($\forall m\in\mathcal{M}$) using either just the sample
paths of $N^{A}(u)$ and $N^{D^{\ast}}(u)$, or just the sample paths of
$N^{B}(u)$ and $N^{D^{\ast}}(u)$. It is clear from
(\ref{three}) that $l_{m}(t)$ may be computed when $\lambda
_{m}(u)$ may
be computed for all $u\in(0,t]$ \textit{and} the sample path of the
counting process for that reaction, $N_{m}(u)$, may be computed over the
same time interval (so that the jump times $\{T_{s}^{m}\leq t\}$ are known).
There are two main elements involved in the argument.

First, the graphical separation $A\bot_{G^{\thicksim}}B|D$ again
has an
important implication for reactants: for any reaction $m$, either $%
R[m]\subseteq A\cup D$ or $R[m]\subseteq B\cup D$.
Recalling~(\ref{cmgm}), only the sample path of the subprocess for the reactants $R[m]$ is
needed to compute $\lambda_{m}(u)$, hence the sample paths of the
subprocesses for either $[A,D]$ or $[B,D]$ suffice, according to
whether $%
R[m]\subseteq A\cup D$ or $R[m]\subseteq B\cup D$. [The sample path of $D$
can clearly be computed from that of $N^{D^{\ast}}(u)$.]

Second, the sample path $(N_{m}(u);u\leq t)$ may be computed using just the
sample paths of $[N^{A},N^{D^{\ast}}]$ or $[N^{B},N^{D^{\ast}}]$,
again according to whether $R[m]\subseteq A\cup D$ or $R[m]\subseteq
B\cup D$%
. To see this, consider each group of reactions in the partition of $%
\mathcal{M}$ given by $[\Delta D_{D},\Delta D_{AB},\Delta(A)\setminus
\Delta(B),\Delta(B)\setminus\Delta(A)]$, beginning with $m\in
\Delta
D_{D}$. By definition the path of $N^{D^{\ast}}(u)$, specifically of
its subcomponent $%
N_{D}^{D}(u)$, allows identification of the jump times corresponding to all
reactions in $\Delta D_{D}$ that change $D$ identically to $m$. But since
such reactions in $\Delta D_{D}$ change $D$ alone, they must do so uniquely
(among reactions in $\Delta D_{D}$) since no 2 columns of $S$ are
equal (Definition
\ref{SKM}). Therefore, the path of $N^{D^{\ast}}(u)$ suffices in this case to
compute $(N_{m}(u);u\leq t)$.

The argument for other groups in the partition
is similar. For $m\in\Delta D_{AB}$, it has already been noted that the
reactions in $\Delta D_{AB}$ change $D$ uniquely (among themselves). The
argument for the last 2 groups is essentially the same. The third
group is
further partitioned as $[\Delta D_{A},\Delta^{\ast}(A)]$, where
$\Delta
^{\ast}(A)$ are the reactions that change $A$ alone. Consider $m\in
\Delta
D_{A}$---again, by definition, the path of $N^{D^{\ast}}(u)$
[specifically, of $%
N_{A}^{D}(u)$] allows identification of the jump times corresponding to the
subset of reactions in $\Delta D_{A}$ that change $D$ identically to $m$.
This subset may now contain more than 1 reaction, but inspection of the
value of the jumps in the sample path of the subprocess for $[A\cup D]$
corresponding to the jump times so identified allows one to ``isolate'' just
those caused by reaction $m$ (since, again, reactions in $\Delta D_{A}$
change $A\cup D$ uniquely among themselves). The argument for $m\in
\Delta
^{\ast}(A)$ is similar, after noting that the jump times of all reactions
in $\Delta^{\ast}(A)$ can be identified by eliminating all those of $%
N_{A}^{D}$ and of $N_{AB}^{D}$.

The preceding two theorems allow the use of Lemma \ref{CIDom} to obtain the
following corollary, which
summarizes the main results of Section \ref{sec4}.
\begin{corollary}
\label{Main3}Let $G$ be the KIG of a standard SKM, $[N,S,\mathsf{P}]$, and
let $[A,B,D]$ be a partition of $\mathcal{V}$. Suppose also that Condition
\ref{IdentGamma} holds for $\Gamma=\Delta(A)\cap\Delta(B)$ (where $%
\Gamma$ is possibly empty). Then
the separation $A\bot_{G^{\thicksim}}B|D$ in the undirected KIG implies
that the global conditional independence $\mathcal{F}_{t}^{A}\indep
{\mathcal{F}_{t}^{B}}|\mathcal{F}_{t}^{D^{\ast}};\mathsf{%
P}_{t}$ holds $\forall t\geq0$, where $\{\mathcal{F}_{t}^{D^{\ast}}\}
$ is
the natural filtration of $N^{D^{\ast}}(t)$.
\end{corollary}
\begin{pf}
Apply Lemma \ref{CIDom} to the 3 $\sigma$-fields $\mathcal{F}_{t}^{A},%
\mathcal{F}_{t}^{B},\mathcal{F}_{t}^{D^{\ast}}\subseteq\mathcal
{F}_{t}^{N}$%
, recalling from Lemma \ref{Likl} that $\mathsf{P}_{t}\ll\mathsf{\tilde
{P}}%
_{t}$. Since $A\bot_{G^{\thicksim}}B|D$, Theorem \ref{Main1} implies
that $\mathcal{F}_{t}^{A}\indep{\mathcal{F}%
_{t}^{B}}|\mathcal{F}_{t}^{D^{\ast}};\mathsf{\tilde{P}}_{t}$. Now
$\mathcal{F}_{t}^{A}\vee\mathcal{F}_{t}^{B}\vee
\mathcal{F}_{t}^{D^{\ast}}=\mathcal{F}_{t}^{N}$, whence $(\mathsf
{dP}_{t}/\break
\mathsf{d\tilde{P}}_{t})|_{\mathcal{F}_{t}^{A}\vee\mathcal
{F}_{t}^{B}\vee
\mathcal{F}_{t}^{D^{\ast}}}=\mathcal{L}_{t}$, which is given by
(\ref{RNderiv}).
Again since $A\bot_{G^{\thicksim}}B|D$, Theorem \ref{Main2}
implies that $\mathcal{L}_{t}=\psi_{AD^{\ast},t}\cdot\psi_{BD^{\ast
},t}$%
, where $\psi_{iD^{\ast},t}$ is a nonnegative, $\mathcal
{F}_{t}^{i}\vee
\mathcal{F}_{t}^{D^{\ast}} $-measurable random variable for $i\in\{
A,B\}$.
Lemma \ref{CIDom} then implies that $\mathcal{F}_{t}^{A}\indep
{\mathcal{F}_{t}^{B}}|\mathcal{F}_{t}^{D^{\ast}};\mathsf{%
P}_{t}$, as required.
\end{pf}

Under the conditions of Corollary \ref{Main3}, the separation $A\bot
_{G^{%
\thicksim}}B|D$ does \textit{not} imply in general that $\mathcal{F}%
_{t}^{A}\indep{\mathcal{F}_{t}^{B}}|%
\mathcal{F}_{t}^{D};\mathsf{P}_{t}$, where the conditioning is now on $%
\mathcal{F}_{t}^{D}$ rather than $\mathcal{F}_{t}^{D^{\ast}}$.
Similarly,\vspace*{1pt}
the separation in the moral graph, $A\bot_{G^{m}}B|D$, does not imply
that $%
\mathcal{F}_{t}^{A}\indep{\mathcal{F}%
_{t}^{B}}|\mathcal{F}_{t}^{D};\mathsf{P}_{t}$. The following theorem and
proof establishes both points. The procedure for constructing $G^{m}$
is the usual one---edges are inserted in the KIG whenever 2 parent
nodes of a common child
are ``unmarried'' (i.e., have no edge between them) and then the
undirected version of the resulting graph is formed.
\begin{theorem}
\label{Countereg}Let $G$ be the KIG of a standard SKM, $[N,S,\mathsf
{P}]$,
and let $[A,B,D]$ be a partition of $\mathcal{V}$. Suppose also that
Condition \ref{IdentGamma}\vspace*{1pt} holds for $\Gamma=\Delta(A)\cap\Delta
(B)$ and
$A\bot_{G^{\thicksim}}B|D$. Then it is possible that neither $\mathcal
{F}%
_{t}^{A}\indep{\mathcal{F}_{t}^{B}}|\mathcal{%
F}_{t}^{D};\mathsf{\tilde{P}}_{t}$ nor $\mathcal{F}_{t}^{A}\indep
{\mathcal{F}_{t}^{B}}|\mathcal{F}_{t}^{D};%
\mathsf{P}_{t}$ holds, where $\{\mathcal{F}_{t}^{D}\}$ is as usual the
internal history of the subprocess $N^{D}(t)$.
\end{theorem}
\begin{pf}
The proof is by example. Consider the standard SKM with $\mathcal{V}%
=\{A,B,D\}$ and reactions
\[
A\rightarrow^{f}D,\qquad D\rightarrow^{r}A,\qquad
D\rightarrow^{\mathit{irr}}B,
\]
which has the KIG, $G=A\longleftrightarrow D\rightarrow B$. Note that
$G^{%
\thicksim}=G^{m}$. Clearly, $\Gamma=\varnothing$ and $A\bot
_{G^{\thicksim}}B|D$. Note also that
$N^{A}(t)=[N_{f}(t),N_{r}(t)]^{\prime}$, $%
N^{D}(t)=[N_{f}(t),N_{r}(t)+N_{\mathit{irr}}(t)]^{\prime}$ and $%
X^{B}(t)-X^{B}(0)=N_{\mathit{irr}}(t)$. It suffices to show that, under both
$\mathsf{%
\tilde{P}}_{t}$ and $\mathsf{P}_{t}$, $\mathsf
{E}[X^{B}(t)-X^{B}(0)|\mathcal{%
F}_{t}^{D}]$ is not a version of $\mathsf{E}[X^{B}(t)-X^{B}(0)|\mathcal
{F}%
_{t}^{A}\vee\mathcal{F}_{t}^{D}]$. First, show that $\mathcal
{F}_{t}^{A}\vee
\mathcal{F}_{t}^{D}=\mathcal{F}_{t}^{N}$. Clearly, $\mathcal
{F}_{t}^{A}\vee
\mathcal{F}_{t}^{D}\subseteq\mathcal{F}_{t}^{N}$, and since $%
N_{\mathit{irr}}(s)=[N_{r}(s)+N_{\mathit{irr}}(s)]-N_{r}(s)$, $N_{\mathit{irr}}(s)$ is measurable $
\mathcal{F}_{t}^{A}\vee\mathcal{F}_{t}^{D}$, hence $\mathcal{F}%
_{t}^{N}\subseteq\mathcal{F}_{t}^{A}\vee\mathcal{F}_{t}^{D}$. It follows
that, under both $\mathsf{\tilde{P}}_{t}$ and $\mathsf{P}_{t}$, $\mathsf
{E}%
[X^{B}(t)-X^{B}(0)|\mathcal{F}_{t}^{A}\vee\mathcal{F}_{t}^{D}]=N_{\mathit{irr}}(t)$
since $N_{\mathit{irr}}(t)$ is measurable $\mathcal{F}_{t}^{A}\vee\mathcal{F}_{t}^{D}
$. However, $N_{\mathit{irr}}(t)$ is clearly not measurable $\mathcal{F%
}_{t}^{D}$ and so cannot be a version of $\mathsf{E}[X^{B}(t)-X^{B}(0)|%
\mathcal{F}_{t}^{D}]$ under either probability measure. In fact, it is
possible to show that under $\mathsf{\tilde{P}}_{t}$, $\mathsf{E}%
[X^{B}(t)-X^{B}(0)|\mathcal{F}_{t}^{D}]=\frac{1}{2}[N_{r}(t)+N_{\mathit{irr}}(t)]$.
\end{pf}

\subsection{Histories of the separator, $D$}
\label{sec43}

It is of interest in applications to understand, for a given partition $
[A,B,D]$ of $\mathcal{V}$, how the histories $\{\mathcal
{F}_{t}^{D^{\ast}}\}
$ and $\{\mathcal{F}_{t}^{D}\}$ differ. A comparison of $N^{D^{\ast}}(t)$
and $N^{D}(t)$ is equivalent to a comparison of the corresponding partitions
of the reactions $\Delta(D)$.
\begin{proposition}
\label{ABDStar}The partition given by $\mathcal{M}^{\ast}(\Delta
(D)):=\{%
\mathcal{M}(\Delta D_{A})\cup\mathcal{M}(\Delta D_{AB})\cup\mathcal{M}
(\Delta D_{B})\cup\mathcal{M}(\Delta D_{D})\}$ is a refinement of the
partition $\mathcal{M}(\Delta(D))$, so that every element of $\mathcal
{M}%
(\Delta(D))$ is a union\vspace*{1pt} of elements of $\mathcal{M}^{\ast}(\Delta(D))$
(see Definition \ref{NADeltaA} for the partition notation used). Hence,
$%
\mathcal{F}_{t}^{D^{\ast}}\supseteq\mathcal{F}_{t}^{D}$ and $\mathcal
{F}%
_{t}^{A}\vee\mathcal{F}_{t}^{B}\vee\mathcal{F}_{t}^{D^{\ast
}}=\mathcal{F}%
_{t}^{N}$ $\forall t$.
\end{proposition}
\begin{pf}
Take an element of $\mathcal{M}(\Delta(D)),\mathcal{M}_{e}(\Delta(D))$
say. Let $m\in\mathcal{M}_{e}(\Delta(D))$ and denote the element of
$%
\mathcal{M}^{\ast}(\Delta(D))$ to which $m$ belongs as $\mathcal{M}%
_{m}^{\ast}(\Delta(D))$. Now $\mathcal{M}_{m}^{\ast}(\Delta
(D))\subseteq
\mathcal{M}_{e}(\Delta(D))$ since all elements of $\mathcal
{M}_{m}^{\ast
}(\Delta(D))$ change $D$ equivalently (resulting in the same change to $D$
as $m$ does). Thus, $\mathcal{M}_{e}(\Delta(D))=\bigcup_{m\in
\mathcal{M}%
_{e}(\Delta(D))}\mathcal{M}_{m}^{\ast}(\Delta(D))$, which
establishes the
first claim. It then follows from Definition \ref{NADeltaA} that
$\mathcal{F}%
_{t}^{D^{\ast}}\supseteq\mathcal{F}_{t}^{D}$ since
elements of $%
N^{D}(t)$ are obtained by summing (where necessary) the appropriate elements
of $N^{D^{\ast}}(t)$. Lemma \ref{Filts} established that $\mathcal{F}%
_{t}^{A}\vee\mathcal{F}_{t}^{B}\vee\mathcal{F}_{t}^{D}=\mathcal{F}_{t}^{N}.
$ But $\mathcal{F}_{t}^{D}\subseteq\mathcal{F}_{t}^{D^{\ast}}$ then
implies $\mathcal{F}_{t}^{N}\subseteq\mathcal{F}_{t}^{A}\vee\mathcal
{F}%
_{t}^{B}\vee\mathcal{F}_{t}^{D^{\ast}}\subseteq\mathcal{F}_{t}^{N}$.
\end{pf}

In computational work with SKMs, establishing if the partitions $%
\mathcal{M}^{\ast}(\Delta(D))$ and $\mathcal{M}(\Delta(D))$ are identical
provides a straightforward means of checking whether the processes $%
N^{D^{\ast}}(t)$ and $N^{D}(t)$ are identical. The two partitions are
identical if and only if there do not exist two reactions in different
elements of $[\Delta D_{A},\Delta D_{AB},\Delta D_{B},\Delta D_{D}]$ that
result in the same change in $D$---that is, there do not exist 2 reactions
in $\Delta(D)$ that change $D$ identically but do not have the same
membership of
both of the sets $[\Delta(A),\Delta(B)]$.
\begin{proposition}
\label{check}Let $[A,B,D]$ be a partition of $\mathcal{V}$, the
species set of an SKM. Then $N^{D^{\ast}}(t)=N^{D}(t)$ $\forall
t,\forall
\omega\in\Omega$, if and only if the following condition holds:
for any 2 reactions $m,\tilde{m}\in\Delta(D)$ with
$S_{m}^{D}=S_{\tilde{m}%
}^{D}$, the reaction $m$ has the same membership of the two
sets $[\Delta(A),\Delta(B)]$ as does the reaction $\tilde{m}$.

Under this condition, $\{\mathcal{F}_{t}^{D^{\ast}}\}=\{\mathcal{F}%
_{t}^{D}\}$.
\end{proposition}
\begin{pf}
If the condition holds both $m$ and $\tilde{m}$ are members of an
equivalence class of
some $\Delta D_{\bullet}$ Hence, any 2 members of an equivalence class
of $%
\Delta(D)$---that is, of an element of $\mathcal{M}(\Delta(D))$---are also
both members of an element of $\mathcal{M}^{\ast}(\Delta(D))$. Therefore,
by Proposition \ref{ABDStar}, $\mathcal{M}(\Delta(D))=\mathcal{M}^{\ast
}(\Delta(D))$, whence $N^{D^{\ast}}(t)=N^{D}(t)$ $\forall t,\forall
\omega
$.

Conversely, suppose $N^{D^{\ast}}(t)=N^{D}(t)$ $\forall t,\forall
\omega$.
Then the vectors $N^{D^{\ast}}(t)=N^{D}(t)$ have the same dimension and
so $\mathcal{M}^{\ast}(\Delta(D))$ cannot be a strict refinement of $%
\mathcal{M}(\Delta(D))$. Hence, by Proposition \ref{ABDStar}, $\mathcal
{%
M}(\Delta(D))=\mathcal{M}^{\ast}(\Delta(D))$. Suppose the
reactions $(m,%
\tilde{m})$ differ in their membership of the two sets $\Delta
(A),\Delta
(B)$. Then $(m,\tilde{m})$ are in different elements of $\mathcal
{M}^{\ast
}(\Delta(D))$ but the same element of $\mathcal{M}(\Delta(D))$, which
is a
contradiction.
\end{pf}

In applications where it is more convenient or practical to include only
internal histories, checking the condition of Proposition
\ref{check}---or equivalently, the equality of $\mathcal{M}(\Delta(D))$ and
$\mathcal
{M}%
^{\ast}(\Delta(D))$---often reveals that the processes $N^{D^{\ast}}(t)$
and $N^{D}(t)$ are similar or identical. This is in part because,
in practice, many elements of $\mathcal{M}(\Delta(D))$ are single
reactions---that is, many of the reactions that change $D$ are uniquely
identified by
the corresponding change in $D$. Furthermore, where $\mathcal{F%
}_{t}^{D}\subset\mathcal{F}_{t}^{D^{\ast}}$ (strictly), the partition
$%
[A,B,D]$ can often be altered slightly to make the processes $N^{D^{\ast
}}(t)$ and $N^{D}(t)$ identical. Examples of this are given in
Section~\ref{sec6} in connection with the red blood cell SKM.

\section{Independence and modularity}
\label{sec5}

Rigorous mathematical definition and identification of
modularizations for biochemical reaction networks is recognized as
being a
difficult problem, especially from a dynamic perspective
\cite{Stelling2006}. A prominent approach has been to
construct a graph representing ``interactions'' between species and to
consider different \textit{partitions} of the species between modules,
maximizing an objective function based on the fraction of edges that are
intra-modular relative to the expected fraction in an ``equivalent,''
randomized graph
when the same partition of species is used \cite{Guimera2005,Alon2005}.
From a stochastic process perspective, the graphs used
often do not encode properly the dependence structure of the molecular
network---for example, in contrast to a KIG, metabolic network graphs typically
omit the local dependence between reactants in the same reaction, only
capturing that between reactant and product. The approach is intended to
operationalize the concept that modules function
``near-independently.''
However, the measure of modularity adopted for the objective function is
rather distant from well-defined notions of dynamic (in)dependence between
species. The local and global conditional independence results
developed in
Sections \ref{sec3} and \ref{sec4} make it possible to add content to and
make rigorous what is meant by near-independence of modules, and to
accommodate ``overlapping'' modules with nonempty intersection.

The term modularization is derived from the biological
literature where ``modularity'' has been much discussed. A modularization here
is a hypergraph of the vertex set of the KIG (i.e., a collection of subsets
of species) with the following property---the internal history at time $t$
of each subset (or module) is conditionally independent of the internal
history of all the other modules, given the history of its intersection with
those modules.
\begin{definition}
\label{Modsatn}Let $\mathcal{V}$ be the species set of an SKM $[N,S,%
\mathsf{P}]$. The finite collection of subsets of $\mathcal{V}$, $%
\{M_{d}|M_{d}\subseteq\mathcal{V}\}$, is a \textit{modularization of
the SKM} if and only if $\bigcup_{d}M_{d}=\mathcal{V}$ and
%
%
\begin{equation} \label{mod}
\mathcal{F}_{t}^{M_{d}}\indep{\mathcal{F}%
_{t}^{\bigcup_{e\neq
d}M_{e}}|\mathcal{F}_{t}^{S_{d}^{\ast}}};\mathsf{P}\qquad
\forall d,t,
\end{equation}
where $S_{d}=M_{d}\cap\{\bigcup_{e\neq d}M_{e}\}$ and the history $\{{%
\mathcal{F}_{t}^{S_{d}^{\ast}}}\}$ is the natural filtration of $N^{{S}
_{d}^{\ast}}(t)$. The latter is given as usual by Definition \ref{NDStar},
applied to the partition $[M_{d}\setminus S_{d},\mathcal{V}\setminus
M_{d},S_{d}]$.
\end{definition}

Note that since $\mathcal{V}\setminus M_{d}=\{\bigcup_{e\neq
d}M_{e}\}\setminus S_{d}$, (\ref{mod}) is equivalent to the
statement $%
\mathcal{F}_{t}^{M_{d}\setminus S_{d}}\indep{\mathcal{F}_{t}^{\mathcal
{V}\setminus M_{d}}|\mathcal{F}%
_{t}^{S_{d}^{\ast}}};\mathsf{P}$ $\forall d,t$. Roughly speaking, the
global evolution on $[0,t]$ of the species in $M_{d}\setminus S_{d}$
and the species
in $\mathcal{V}\setminus M_{d}$ (``the rest of the network'') are
conditionally independent given the history of the intersection,
${\mathcal{F%
}_{t}^{S_{d}^{\ast}}}$. We will say that two modularizations are
\textit{nested} if each module of one of the modularizations is contained in some
module of the other modularization.

Of course some modularizations of an SKM will be more useful than
others. It
will usually be desirable for the intersections $S_{d}$ to contain a
relatively small number of species and to be able to move between nested
modularizations, thus considering finer and coarser levels of resolution.
Computationally efficient methods are developed below for the identification
of such modularizations that are based around the maximal prime
decomposition of the undirected version of the KIG of the SKM,
$G^{\thicksim}$. It will be proved below that applying such graphical
decomposition
methods results in subgraphs whose vertex sets, $\{M_{d}\}$ say,
satisfy the
graphical separation $M_{d}\perp_{G^{\thicksim}}\bigcup_{e\neq
d}M_{e}{|}%
S_{d}$ $\forall d$. Therefore, under the
conditions of Corollary \ref{Main3},
the required global dynamic independence of (\ref{mod})
holds for all $%
d$, and $\{M_{d}\}$ constitutes a modularization according to
Definition \ref{Modsatn}.

\subsection{Identifying modularizations by graph decomposition}\label{sec51}

Some definitions from the graphical literature will prove useful (for
further details,
see \cite{Lauritzen96}). An
undirected graph is said to be \textit{complete} if there is an edge
between all pairs of vertices in its vertex set. Let $H$ be an
undirected graph
with vertex set $\mathcal{V}$. The subgraph
induced by $M_{d}\subset\mathcal{V}$, $H(M_{d})$, consists of the vertices
in $M_{d}$ and exactly the edges between those vertices that occur in $H$
itself. A partition $[A,B,D]$ of $\mathcal{V}$, $A,B\neq\varnothing$, forms
a \textit{decomposition} of $H$ into the subgraphs $H(A\cup D)$ and
$H(B\cup
D)$ if the separation $A\bot_{H}B|D$ holds and the subgraph $H(D)$ is
complete. The subgraph $H(M_{d})$ is \textit{prime} if there does not exist
a decomposition of $H(M_{d})$.
\begin{definition}
\label{MPD} Let $H$ be an undirected graph with vertex set $\mathcal{V}$,
and \mbox{$M_{d}\subseteq\mathcal{V}$}. The induced subgraph $H(M_{d})$ is a
maximal prime subgraph of $H$ if $H(M_{d})$ is prime and there exists a
decomposition of $H(N)$ for all $N$ satisfying $M_{d}\subset N\subseteq
\mathcal{V}$. The \textit{maximal prime subgraph decomposition}
(\textit{MPD}) of $H$ is given by $\{H(M_{d})\}$, the unique collection of maximal prime
subgraphs of $H$, and satisfies that $\bigcup_{d}M_{d}=\mathcal{V}$.
\end{definition}

A \textit{junction tree} representation of the MPD,
$\mathcal{T}_{\mathrm{MPD}}$, always exists and has the subsets $\{M_{d}\}$ as its
clusters (i.e., as the vertices of the junction tree)
\cite{Olesen2002}. A junction tree $\mathcal{T}$ is a connected, undirected
graph without cycles in which the intersection
of any 2 clusters of the tree, $M_{d}\cap M_{e}$ $(d\neq e)$, is contained
in every cluster on the unique path in $\mathcal{T}$ between $M_{d}$
and $M_{e}$. Such trees will prove very useful in visualizing, representing
and manipulating modularizations of SKMs. We say, for reasons that will become
apparent, that any 2 clusters adjacent in the tree are separated
by their intersection, and call that intersection a \textit{separator}
\textit{of} $\mathcal{T}$.

The SKM modularization algorithm presented below contains as a special case
the method due to \cite{Olesen2002} for computation of $\mathcal
{T}%
_{\mathrm{MPD}}$, applied to the undirected version of the KIG, $G^{\thicksim}$.
The advantage of this version of Algorithm \ref{Tmod} is that it can be
fully automated to identify the MPD modularization of the SKM in a manner
that is computationally feasible even for very large SKMs. However, it will
often be informative to consider a range of nested modularizations in order
to explore the different levels of organization of the reaction
network. To
this end, the general version of Algorithm \ref{Tmod} first obtains a
junction tree of the clique decomposition for $G_{T}^{\thicksim}$ (a
minimal triangulation of $G^{\thicksim})$---this provides the finest, most
detailed modularization that is identified. The clique decomposition of
$G_{T}^{\thicksim}$ is unique (since it corresponds to the MPD of
$G_{T}^{\thicksim}$). Coarser-grained modularizations,
including the MPD one, are obtained by successively aggregating adjacent
clusters in the junction tree.
\begin{algorithm}
\label{Tmod}Let $G$ be the KIG of an SKM.

\begin{enumerate}
\item Construct $G^{\thicksim}$, the undirected version of $G$;

\item Construct $G_{T}^{\thicksim}$, a minimal triangulation of $G^{%
\thicksim}$;

\item Obtain the clique decomposition of $G_{T}^{\thicksim}$ with the
cliques, $\{C_{1},C_{2},\ldots,C_{\delta}\}$ say, ordered to satisfy the
running intersection property (i.e., for $e=2,\ldots,\delta,\exists
d^{\ast
}\in\{1,\ldots,e-1\}$ s.t. $C_{e}\cap\{\bigcup_{i=1}^{e-1}C_{i}\}\subseteq
C_{d^{\ast}}$);

\item Organize the clique decomposition as a (rooted) junction tree
$\mathcal{T}_{C}$ in which, for $e=2,\ldots,\delta$, the parent of
$C_{e}$ is
$C_{d^{\ast}};$ set $\mathcal{T}=\mathcal{T}_{C}$;

\item Either go to step 7 or,
select a pair of adjacent clusters $(C_{i},C_{j})$ in $\mathcal
{T}$ $(i<j)$ and \textit{aggregate} them by updating
$\mathcal{T}$ as follows:
set $P=\mathrm{pa}(C_{i})$ and $C=\{\operatorname{ch}(C_{i})\cup \operatorname{ch}(C_{j})\}\setminus C_{j}$,
replace cluster $i$ by $C_{i}\cup C_{j}$ (retaining its numbering,
$i$), set
$\mathrm{pa}(C_{i})=P$, and set $\operatorname{ch}(C_{i})=C;$

\item Go to step 5;

\item Return $\mathcal{T}_{\mathrm{MOD}}=\mathcal{T}$.
\end{enumerate}
\end{algorithm}

The property that $G_{T}^{\thicksim}$ is triangulated is
equivalent to saying that $G_{T}^{\thicksim}$ can be decomposed
recursively until all the resulting subgraphs are complete \cite{Lauritzen96}.
Such a recursive decomposition produces a collection of subgraphs
containing the cliques $\{G_{T}^{\thicksim}(C_{d})\}$, that is, the maximally
complete subgraphs of $G_{T}^{\thicksim}$. Triangulation refers to the
operation of adding
edges to $G^{\thicksim}$ so that it becomes triangulated. The
triangulation $G_{T}^{\thicksim}$
in step 2 must be minimal---that is, one for which removal of any edge added
during triangulation results in an untriangulated graph---otherwise,
Remark~\ref{rem51} below does not
hold, in general.

Efficient algorithms have been developed in the graphical literature
for both minimal
triangulation and clique decomposition (see
\cite{Olesen2002,CDLS2007}) which can be exploited here to compute the
SKM modularizations
and associated junction trees. The following special case of Algorithm
\ref{Tmod} returns the junction tree representation of the maximal prime
decomposition (MPD) of the undirected KIG, $G^{\thicksim}$
\cite{Olesen2002}.
\begin{remark}\label{rem51}
Algorithm \ref{Tmod} returns $\mathcal{T}_{\mathrm{MPD}}$ for the undirected
KIG, $G^{\thicksim}$, when step 5 is replaced by:

5$^{\prime}$. While [there exists a separator $S$ of
$\mathcal{T}$ such that $G^{\thicksim}(S)$ is incomplete],
aggregate within $\mathcal{T}$ the 2 clusters separated by $S$;
then go to step 7.
\end{remark}

It is worth noting the time complexity of steps 2 and 4. The general
problem of finding an optimal triangulation of an undirected graph
(i.e., one that adds least edges among all triangulations) is
\textit{NP-hard}. The complexity of minimal triangulation (step 2) is
$\mathcal{O}(ne)$ where $e$ is the number of edges in $G^{\thicksim}$,
\cite{Olesen2002}. The complexity of constructing the clique junction
tree $\mathcal{T}_{C}$ (steps 3 and 4 combined) is $\mathcal
{O}(n^{2})$, \cite{Olesen2002}.

\subsection{Nested modularizations and junction trees}\label{sec52}

A concise proof that the clusters of the tree $\mathcal
{T}_{\mathrm{MOD}}$
returned by Algorithm \ref{Tmod} constitute a modularization of the
SKM---with any choice of aggregation scheme in stage 5---is made
possible by
establishing that $\mathcal{T}_{\mathrm{MOD}}$, like $\mathcal
{T}_{C}$, is a junction tree,
and that the intersections of adjacent clusters of $\mathcal
{T}_{\mathrm{MOD}}$ continue to
correspond to separators in $G_{T}^{\thicksim}$, and hence in $%
G^{\thicksim}$. The following proposition does just that.
\begin{proposition}
\label{JTProp}Let $\mathcal{T}_{\mathrm{MOD}}$ be the undirected graph
returned by applying Algorithm \ref{Tmod} to the KIG, $G$, of an SKM.
Denote the clusters (modules) of $\mathcal{T}_{\mathrm{MOD}}$ by $%
\{M_{d}\}$. Then $\mathcal{T}_{\mathrm{MOD}}$ is a junction tree. Suppose
that $(M_{d},M_{e})$ are any 2 adjacent clusters in $\mathcal{T}%
_{\mathrm{MOD}}$ with separator $S_{\mathit{de}}:=M_{d}\cap M_{e}$, and that (as is
conventional) the edges $M_{d}\sim M_{e}$ are labeled by the corresponding
separator $S_{\mathit{de}}$.

Then $S_{\mathit{de}}=V_{\mathit{de}}\cap V_{\mathit{ed}}$ and the graphical separation
$V_{\mathit{de}}\perp
V_{\mathit{ed}}|S_{\mathit{de}}$ holds in $G_{T}^{\thicksim}$, and hence in $G^{\thicksim}$,
where $V_{\mathit{de}}$ $(V_{\mathit{ed}})$ is the union of the clusters in $\mathcal{T%
}_{\mathrm{MOD}}^{\mathit{de}}$ ($\mathcal{T}_{\mathrm{MOD}}^{\mathit{ed}}$), the $\mathcal{T}%
_{\mathrm{MOD}}^{\bullet}$ are the 2 subtrees obtained by cutting the edge $%
M_{d}\sim M_{e}$ in $\mathcal{T}_{\mathrm{MOD}}$, and $M_{d}\subseteq V_{\mathit{de}}$
($M_{e}\subseteq V_{\mathit{ed}}$).
\end{proposition}

Proof of Proposition \ref{JTProp} is given in Appendix \ref{Prooffs}.

We can now state and prove the result that establishes the validity of our
modularization identification methods.
\begin{theorem}
\label{JTMod} Let $G$ be the KIG of a standard SKM, $[N,S,\mathsf{P}]$, and
let $\mathcal{T}_{\mathrm{MOD}}$ be the junction tree of modules, $\{M_{d}\}$,
returned by Algorithm \ref{Tmod}. Suppose
also that Condition \ref{IdentGamma} holds for $\Gamma_{d}=\Delta
(M_{d}\setminus S_{d})\cap\Delta(\mathcal{V}\setminus M_{d})$
$\forall d$. \textit{Then} $\{M_{d}\}$ is a modularization of the SKM
in the
sense of
Definition \ref{Modsatn}; and each $S_{d}$ is given by $\bigcup_{e\in
ne(M_{d})}S_{\mathit{de}}$ [where $ne(M_{d})$ is the indices of those
clusters that have edges with $M_{d}$ in $\mathcal{T}_{\mathrm{MOD}}$].
Furthermore, each module residual $M_{d}\setminus S_{d}$ is locally
independent of
$\mathcal{V}\setminus M_{d}$ given the internal history of $M_{d}$.
\end{theorem}
\begin{pf}
By Corollary \ref{Main3}, it suffices to show
that the separation $\{M_{d}\setminus S_{d}\} \perp
_{G_{T}^{\thicksim}}\{%
\mathcal{V}\setminus M_{d}\}|S_{d}$ holds in $G_{T}^{\thicksim}$, for
all $d$, since then $\{M_{d}\setminus S_{d}\}\perp_{G^{\thicksim}}\{%
\mathcal{V}\setminus M_{d}\}|S_{d}$ holds in the undirected KIG
$G^{\thicksim}$.
This follows because every path in $G^{\thicksim}$ from $%
M_{d}\setminus S_{d}$ to $\mathcal{V}\setminus M_{d}$ is also such a path
in $G_{T}^{\thicksim}$. Recall that by definition $S_{\mathit{de}}:=M_{d}\cap M_{e}.
$ Hence,
%
%
\begin{eqnarray}
\label{sde}
S_{d} &=&\biggl\{\bigcup_{e\in ne(M_{d})}(M_{d}\cap M_{e})\biggr\}\cup
\biggl\{\bigcup_{e\notin ne(M_{d})}(M_{d}\cap M_{e})\biggr\}
\nonumber\\[-8pt]\\[-8pt]
&=&\bigcup_{e\in ne(M_{d})}(M_{d}\cap M_{e})=\bigcup_{e\in ne(M_{d})}S_{\mathit{de}},
\nonumber
\end{eqnarray}
where the second line holds by the fact that $\mathcal
{T}_{\mathrm{MOD}}$ is
a junction tree (Proposition~\ref{JTProp}) since, for $e\notin
ne(M_{d})$, $(M_{d}\cap M_{e})$ is contained in $M_{\tilde{e}}$, and
thus in $%
S_{d\tilde{e}}$ for some $\tilde{e}\in ne(M_{d})$ lying on the unique path
between $M_{d}$ and $M_{e}$ in $\mathcal{T}_{\mathrm{MOD}}$.

By Proposition \ref{JTProp}, $M_{d}\perp_{G_{T}^{\thicksim%
}}V_{\mathit{ed}}|S_{\mathit{de}}$ $\forall e\in ne(M_{d})$, since $M_{d}\subseteq V_{\mathit{de}}$.
Hence, $M_{d}\perp_{G_{T}^{\thicksim}}V_{\mathit{ed}}|\{\bigcup_{e\in
ne(M_{d})}S_{\mathit{de}}\}$ and, since this holds for all $e\in ne(M_{d})$, $%
M_{d}\perp_{G_{T}^{\thicksim}}\{\bigcup_{e\in
ne(M_{d})}V_{\mathit{ed}}\}|\{\bigcup_{e\in ne(M_{d})}S_{\mathit{de}}\}$. Now $\{\bigcup
_{e\in
ne(M_{d})}V_{\mathit{ed}}\}=\{\bigcup_{e\neq d}M_{e}\}$. To see this, note that the
latter is the union of those clusters reachable by paths in $%
\mathcal{T}_{\mathrm{MOD}}$ that start with the edge $M_{d}\sim
M_{e}$ for
some $e\in ne(M_{d})$ (since $\mathcal{T}_{\mathrm{MOD}}$ is
connected); and
$V_{\mathit{ed}}$ is the union of those clusters reachable by paths in
$\mathcal{T}_{\mathrm{MOD}}$ that start at the node $M_{e}$ (since
$\mathcal{T}_{\mathrm{MOD}}^{\mathit{ed}}$ is connected). Therefore,
using (\ref{sde}), $%
M_{d}\perp_{G_{T}^{\thicksim}}\{\bigcup_{e\neq d}M_{e}\}|S_{d}$, as
required.
\end{pf}

\section{SKM of red blood cell}
\label{sec6}

We now apply the modularization techniques of Section \ref{sec5}
and the
underlying dynamic independence theory on which they are based to identify
biologically interesting modulariszations of an SKM of the human red blood
cell. The study of this metabolic reaction network was an
early success of a systems biology approach
\cite{Schuster95,Palsson2002}. There now exists detailed knowledge of
the component
reactions as a result of at least three decades of
research on both the biochemical and mathematical modeling fronts. The
identification of aggregates of metabolites (i.e., species) and regulatory
structures in the red blood cell has also received attention from a systems
biology perspective \cite{Palsson2002,Palsson2003}. This particular
reaction network therefore constitutes a suitable test-bed to establish the
utility and applicability of our approach. In contrast to this work,
\cite{Palsson2002} aims to identify ``pools'' of metabolites in the red
blood cell,
that is ``aggregate groups of [species] which [\ldots] move together in a
concerted manner,'' rather than groups that move independently given an
appropriate conditioning set of species.

The SKM studied is the one implied by the metabolic network of the red blood
cell published in the open access Biomodels Database \cite{BIOMOD2006},
which in turn is a slightly extended version of the kinetic model of
\cite{Schuster95} and \cite{Holzhutter2004}. The SKM consists of 38 reactions,
with 45 different biochemical species in the species set $\mathcal{V}$ (the
enzymes, i.e., catalysts, involved are omitted from $\mathcal{V}$ as
they do not appear
explicitly in the reaction mechanisms). Full details are available from
\cite{BIOMOD2006}. The direction of the reactions is as for the kinetic
model in Table 1 of \cite{Holzhutter2004}, except for 8 additional
reactions which are all included as dissociation reactions. It was verified
that the SKM, henceforth $\mathcal{SKM}_{\mathrm{rbc}}$, is a standard SKM (according
to Definition~\ref{stdSKM}). The names of the biochemical species in $
\mathcal{V}$ and the associated abbreviations used are given in
Appendix \ref{species}. For details of the reactions in $\mathcal{M}$, the reader is
referred to~\cite{BIOMOD2006}.

%
%
\begin{figure}

\includegraphics{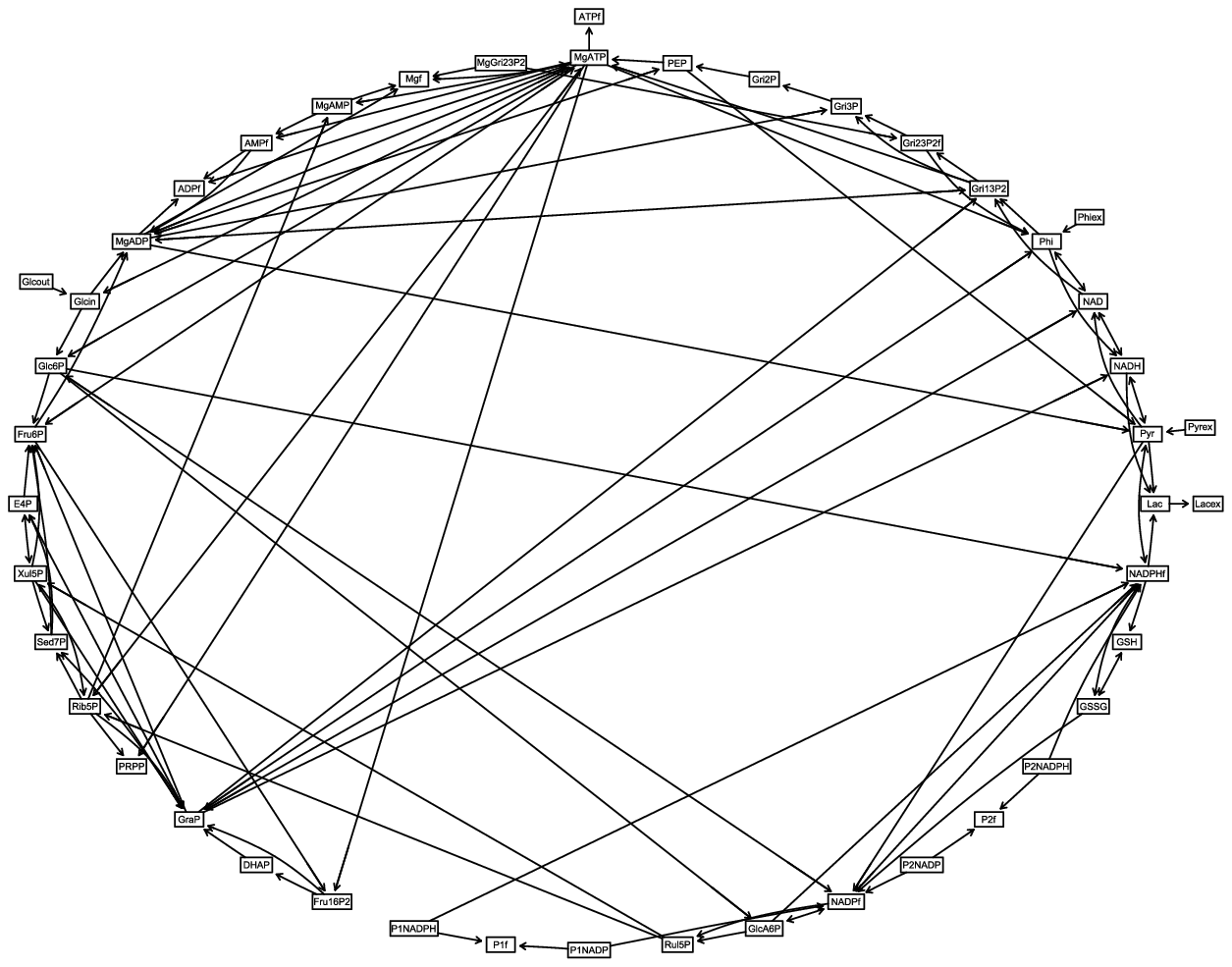}

\caption{Kinetic independence graph (KIG) of
$\mathcal{SKM}_{\mathrm{rbc}}$, the Metabolic Network of the Human Red Blood Cell
\protect\cite{BIOMOD2006}. The KIG is constructed according to
Definition \protect\ref{LIGnew}. Full species names are given in
Appendix \protect\ref{species}.}\label{KIGEry}
\end{figure}

\begin{figure}

\includegraphics{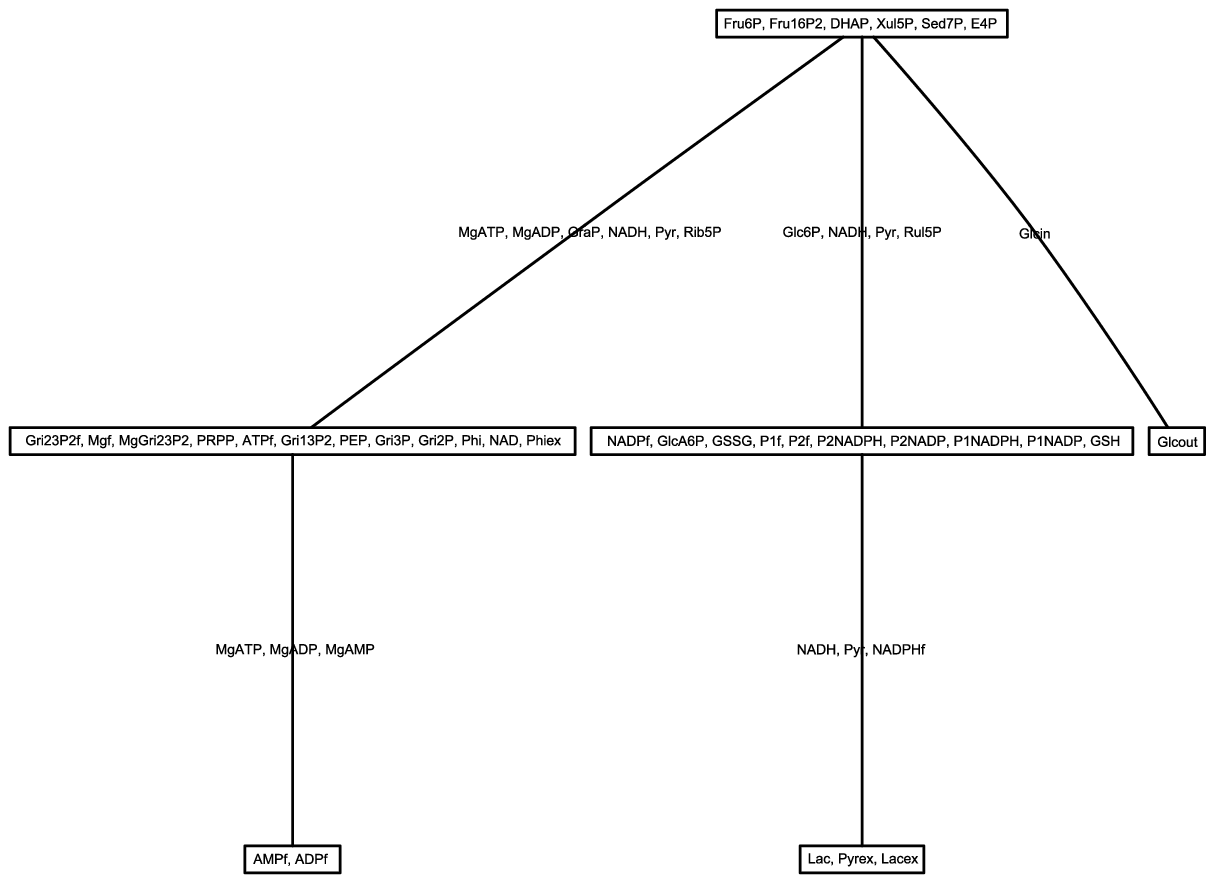}

\caption{Junction tree representation, $\mathcal{T}_{\mathrm{MOD},1}$, of a
modularisation of $\mathcal{SKM}_{\mathrm{rbc}}$, the Metabolic Network of the
Human Red Blood Cell \protect\cite{BIOMOD2006}. The global dynamic
independences $\mathcal{F}_{t}^{M_{d}}\perp\hspace*{-5.2pt}\perp{\mathcal{F}_{t}^{\bigcup
_{e\neq d}M_{e}}|\mathcal{F}_{t}^{S_{d}^{\ast}}}$ hold for each module
$M_{d}$ (see Definition \protect\ref{Modsatn}). The modules
(rectangles) are
labeled with their residuals and edges are labeled with the
intersection of adjacent modules. Full species names are given in
Appendix \protect\ref{species}.}\label{TMod2}
\end{figure}

Figure \ref{KIGEry} depicts the kinetic independence graph $G$ for
$\mathcal{%
SKM}_{\mathrm{rbc}}$. The graph is a powerful visual aid to understanding the
architecture of the molecular network and can be preliminarily
inspected for
interesting local independences and separations in the undirected
version $%
G^{\thicksim}$. The clique decomposition, $\mathcal{T}_{C}$,
from Algorithm \ref{Tmod} for $\mathcal{SKM}_{\mathrm{rbc}}$ has many clusters (20
out of 38) for which $M_{d}\setminus S_{d}$ is the empty set. It is
therefore desirable to implement Algorithm \ref{Tmod} with a substantial
degree of pairwise cluster aggregation in step 5. On the other hand,
$\mathcal{T}_{\mathrm{MPD}}$ for this SKM is overly coarse-grained for most
purposes. Figure \ref{TMod2} depicts a particular junction tree
$\mathcal{T}_{\mathrm{MOD},1}$ returned by Algorithm \ref{Tmod}, with the choice of
aggregations guided both by the structure of $\mathcal{T}_{C}$
itself and the goal of a modularisation that offers biological
insight. This approach relies on and takes advantage of the flexibility
offered by Algorithm \ref{Tmod}---some exploration of alternative
modularizations by the user is required, but no prior information about
possible modularizations is needed.

%

The junction tree $\mathcal{T}_{\mathrm{MOD},1}$ in Figure \ref{TMod2} is
labelled as follows. The $d$th module (rectangle) is labeled with the
species in the ``residual module,'' $M_{d}\setminus S_{d}$, and each
edge, $%
M_{d}\sim M_{e}$, is labeled with the species in the separator $S_{\mathit{de}}$,
that is the intersection of the modules connected by that edge. It was
verified that, for all $d$, Condition \ref{IdentGamma} holds for $\Gamma
_{d}=\Delta(M_{d}\setminus S_{d})\cap\Delta(\mathcal{V}\setminus M_{d})$
, as required by Theorem~\ref{JTMod}. Recall that Proposition \ref{JTProp}
implies $S_{\mathit{de}}=V_{\mathit{de}}\cap V_{\mathit{ed}}$ and $V_{\mathit{de}}\perp_{G^{\thicksim%
}}V_{\mathit{ed}}|S_{\mathit{de}}$ in $G^{\thicksim}$. Such separations may be
conveniently read off from any junction tree $\mathcal
{T}_{\mathrm{MOD}}$ since $S_{\mathit{de}},V_{\mathit{de}}$
and $V_{\mathit{ed}}$ are all immediately apparent from
examination of the tree. Similarly, the defining conditional independencies
of the modularization, namely $\mathcal{F}_{t}^{M_{d}\setminus
S_{d}}\indep
{\mathcal{F}_{t}^{\{\bigcup_{e\neq
d}M_{e}\}\setminus S_{d}}|\mathcal{F}_{t}^{S_{d}^{\ast}}};\mathsf{P}$
$\forall d$ (\ref{mod}), may be read off the junction tree
using Theorem
\ref{JTMod}, $S_{d}$ being given by the union of the labels of all edges
that connect with the $d$th module.

Having obtained a modularization such as $\mathcal
{T}_{\mathrm{MOD},1}$, the
next stage is to ask what are the interesting features that emerge from a
biochemical and systems biological perspective. Each of the main
modules of
$\mathcal{T}_{\mathrm{MOD},1}$ turns out to contain like species,
either in terms of
their molecular structure (e.g., the groupings of
monosaccharide-phosphate sugar molecules and phosphoglycerate
molecules) or their function (e.g., the grouping of species
involved in reduction--oxidation reactions), or both.

Specific modules and their
residuals are denoted by their first constituent species in the subsequent
discussion. Consider first the central residual\break $\{\mathit{NAD
Pf},\ldots\}$ in $\mathcal{T}_{\mathrm{MOD},1}$, the residual of what will be termed the
\textit{Redox module} (for \textit{Red}uction-\textit{ox}idation). The red
blood cell is subject to oxidative stress due to reactive oxygen species,
which if left unchecked leads to cell lysis (bursting) and consequent
anemia. All of this
residual's species can be seen to play a role in the control of such
oxidative stress. Glutathione (\textit{GSH}) acts as an antioxidant,
scavenging reactive oxygen species and itself being oxidised as a result
(giving rise to the reaction $2\mathit{GSH}\rightarrow\mathit{GSSG}$).
The cell must maintain adequate levels of \textit{GSH}, which it does by
producing large amounts of \textit{NADPH} for use in the reduction of
\textit{GSSG} (by the reaction $\mathit{GSSG}+\mathit{
NADPH}\rightarrow2\mathit{GSH}+
\mathit{NADP}$). Production of \textit{NADPH} is via 2 reactions
(usually described as the oxidative phase of
the pentose phosphate pathway), both of
which involve \textit{GlcA6P}. Both \textit{NADP} and \textit{NADPH} are
also found bound to the proteins $P1$ and $P2$. Notice that the reduced
forms \textit{NADPH} and \textit{NADH} are both found in the module's
separator (edge) with $\{\mathit{Lac},\mathit{Pyrex},\mathit{Lacex}\}$, since both influence
the intensity of lactate $(\mathit{Lac})$ production and export as reactants for the
reduction of pyruvate $(\mathit{Pyr})$.

The module $\{\mathit{NADPf},\ldots\}$ clearly has an important function
in oxidative
stress control and in reduction--oxidation reactions more generally
within the red blood cell.
Of course, these functions of its individual species are well known. That
their dynamic evolution, together with that of lactate $(\mathit{Lac})$, is globally
independent of all the other species in the network conditional on the
internal history of $\{$\textit{Fru6P}, \textit{Glc6P}, \textit{NADH},
\textit{Pyr}, \textit{Rul5P}, \textit{Pyrex}$\}$
is an
insight provided by the modularizations (see also the derivation of
$\mathcal{T}_{\mathrm{MOD},2}$ below). Assigning function(s) where possible to each
module of a given modularization, $\mathcal{T}_{\mathrm{MOD}}$, is
likely to
improve both understanding of a network and ultimately aid attempts to
control it. For reasons of space, comments related to
the remaining two large residuals of $\mathcal{T}_{\mathrm{MOD},1}$ may be
found as part of the discussion of $\mathcal{T}_{\mathrm{MOD},2}$ below.

%
%
\begin{figure}[b]

\includegraphics{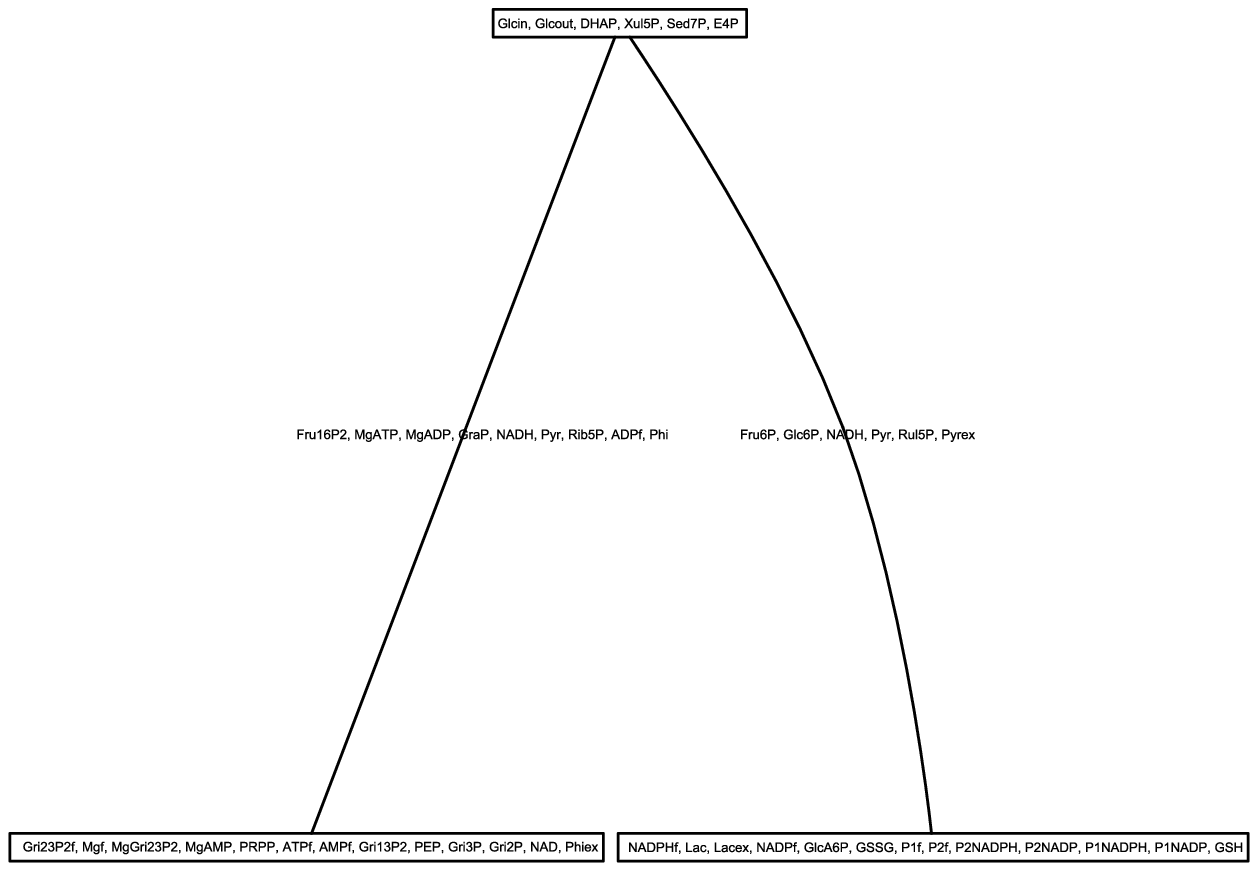}

\caption{Junction tree representation, $\mathcal{T}_{\mathrm{MOD},2}$, a
coarser-grained modularization of $\mathcal{SKM}_{\mathrm{rbc}}$, the Metabolic
Network of the Human Red Blood Cell \protect\cite{BIOMOD2006}. The
global dynamic independences
$\mathcal{F}_{t}^{M_{d}}\perp\hspace*{-5.2pt}\perp{\mathcal{F}_{t}^{\bigcup_{e\neq
d}M_{e}}|\mathcal{F}_{t}^{S_{d}}}$ hold for each module $M_{d}$ (see
Definition \protect\ref{Modsatn}). The rectangles contain module
residuals and
edges are labeled with the intersection of adjacent modules. Full
species names are given in Appendix \protect\ref{species}.}\label{TMod3}
\end{figure}

The structure of $\mathcal{T}_{\mathrm{MOD},1}$ encourages further
aggregation in an obvious manner. A second modularization,
$\mathcal{T}_{\mathrm{MOD},2}$, of $\mathcal{SKM}_{\mathrm{rbc}}$ is thus shown in Figure \ref{TMod3}.
[It was verified that, in this case also, Condition \ref{IdentGamma} holds
for $\Gamma_{d}=\Delta(M_{d}\setminus S_{d})\cap\Delta(\mathcal{V}%
\setminus M_{d})$ $\forall d$.] $\mathcal{T}_{\mathrm{MOD},2}$ may be
derived from $\mathcal{T}_{\mathrm{MOD},1}$in two steps. First,
the modules $\{\mathit{AMPf},\mathit{ADPf}\}$, $\{\mathit{Lac},\mathit{Pyrex},\mathit{Lacex}\}$ and $%
\{\mathit{Glcout}\}$ are aggregated with their adjacent modules in
$\mathcal{T}_{\mathrm{MOD},1}$. Second, a small number of species in residuals are then
judiciously
included also in an additional module so that they both fall instead in the
relevant separator and the condition of Proposition \ref{check} is satisfied
for each partition $[M_{d}\setminus S_{d},\mathcal{V}\setminus
M_{d},S_{d}]$. By Proposition \ref{check}, this ensures that
$N^{S_{d}}$ and
$N^{S_{d}^{\ast}}$ are the same subprocess, whence $\{{\mathcal{F}%
_{t}^{S_{d}}\}=\{\mathcal{F}_{t}^{S_{d}^{\ast}}\}}$ for $d=1,2,3$. Clearly,
the second step need only be performed if it is desired to be able to
replace ${\mathcal{F}_{t}^{S_{d}^{\ast}}}$ by ${\mathcal{F}_{t}^{S_{d}}}$
in the defining conditional independencies of the modularization (Definition
\ref{Modsatn}). The species involved in this case are $\{\mathit{
Fru6P},\mathit{Fru16P2},\mathit{Phi}, \mathit{ADPf},\mathit{Pyrex}\}$ and these therefore now appear in
the edge
labels (separators) of $\mathcal{T}_{\mathrm{MOD},2}$ rather than in the
residuals. The proposition below establishes that the validity of the
modularization remains unchanged by such an operation.
\begin{proposition}
\label{copied}Suppose that $\{M_{d}|M_{d}\subseteq\mathcal{V}\}$, is a
modularization according to Definition \ref{Modsatn} of a standard SKM
$[N,S,%
\mathsf{P}]$, and that the modularization satisfies, for all $d$, the
separation $\{M_{d}\}\perp_{G^{\thicksim}}\{\bigcup_{e\neq d}M_{e}\}|S_{d}$
in the undirected KIG $G^{\thicksim}$. Define a new collection of
subsets $%
\{\tilde{M}_{d}|\tilde{M}_{d}\subseteq\mathcal{V}\}$ where $\tilde{M}%
_{d}=M_{d}\cup\{\bigcup_{e\neq d}c_{\mathit{ed}}\}$ and, $\forall e\neq d$, $%
c_{\mathit{ed}}\subset M_{e}$ and $c_{\mathit{ed}}\cap M_{d}=\varnothing$
$(c_{\mathit{ed}}=\varnothing$
being allowed). The species $c_{\mathit{ed}}$ are called those ``copied from $e$
to $d$.''

Then $\{\tilde{M}_{d}\}\perp_{G^{\thicksim}}\{\bigcup_{e\neq d}\tilde
{M}%
_{e}\}|\tilde{S}_{d}$ $\forall d$ and, provided that
Condition \ref{IdentGamma} continues to hold for $\tilde{\Gamma}_{d}=\Delta(\tilde{M}
_{d}\setminus\tilde{S}_{d})\cap\Delta(\mathcal{V}\setminus\tilde
{M}%
_{d})$ $\forall d$, $\{\tilde{M}_{d}|\tilde{M}_{d}\subseteq\mathcal
{V}\}$
is also a modularization of the SKM $[N,S,\mathsf{P}]$.
\end{proposition}
\begin{pf}
Clearly, $\bigcup_{d}\tilde{M}_{d}=\mathcal{V}$. By Corollary
\ref{Main3}, it
suffices to show that $\{\tilde{M}_{d}\}\perp_{G^{\thicksim}}\{\bigcup
_{e\neq d}\tilde{M}_{e}\}|\tilde{S}_{d}$ $\forall d$. Let $%
t_{d}:=\bigcup_{e\neq d}c_{\mathit{ed}}$, the species copied to $d$, and $%
f_{d}:=\bigcup_{e\neq d}c_{\mathit{de}}$, the species copied from $d$. The
separation $%
\{M_{d}\}\perp_{G^{\thicksim}}\{\bigcup_{e\neq d}M_{e}\}|S_{d}$ implies
that $\{M_{d}\cup t_{d}\}\perp_{G^{\thicksim}}\{\bigcup_{e\neq
d}M_{e}\}\cup f_{d}|\{S_{d}\cup t_{d}\cup f_{d}\}$, which yields the
required result since $\tilde{S}_{d}=\{M_{d}\cup t_{d}\}\cap\lbrack\{
\bigcup
_{e\neq d}M_{e}\}\cup f_{d}]=S_{d}\cup t_{d}\cup f_{d}\cup\varnothing$.
\end{pf}

There are 3 modules comprising $\mathcal{T}_{\mathrm{MOD},2}$ which
together contain 45 different species, of which 32 distinct species are
found only in module
residuals (and hence are found in exactly 1 residual). The redox module
$\{\mathit{NADPHf},\ldots\}$ has already been discussed above. The module $%
\{\mathit{Glcin},\ldots\}$ has the largest intersection with the rest of the
network and
acts as a linking module; it will be termed the \textit{MPS}
(\textit{Monosaccharide-Phosphate Sugar} module). The two modules $\{\mathit{NADPHf}
,\ldots\}$ and $\{\mathit{Gri23P2f},\ldots\}$, by contrast, have only 2
species in common,
namely ($\mathit{NADH},\mathit{Pyr}$)---these are the only species common to all
three modules. The
\textit{MPS} residual contains species that all belong to a single
chemical class of molecule, namely monosaccharide sugar molecules
(mostly with phosphate groups attached), with a further 6 different
monosaccharide-phosphates (\textit{MPs}) found in the rest of the module.
Interestingly, the \textit{MPs} of the module are those found in two
``pathways'' traditionally discussed separately---the pentose phosphate and
glycolytic pathways. Indeed, the \textit{MPs} (\textit{Glc6P}, \textit{Fru6P}, \textit{GraP}) all
participate in reactions found in both ``pathways.''

The third and final module $\{\mathit{Gri23P2f},\ldots\}$ will be termed
the \textit{PGA} (\textit{Phos- phoGlycerate-Adenosine}) module, according to the chemical class
of some of its constituents. It contains all of the
phosphoglycerate molecules in the species set~$\mathcal{V}$, namely
(\textit{Gri23P2f}, \textit{MgGri23P2f}, \textit{Gri13P2},
\textit{Gri3P}, \textit{Gri2P}), together also with
all of the
adenosine phosphate molecules (\textit{ATPf}, \textit{ADPf}, \textit{AMPf})---both free
and complexed with
magnesium (\textit{Mg}). The module also contains all of the species
involved in reactions of the
so-called ``pay-off phase'' of glycolysis whose function is the
production of
the high-energy compounds \textit{ATP} and \textit{NADH}. That the
dynamic evolution of,
for example, all phosphoglycerates together with \textit{PEP} is
globally independent of
all the other species in $\mathcal{V}$ conditional on the internal
history of (\textit{Fru16P2}, \textit{MgATP}, \textit{MgADP},
\textit{GraP}, \textit{NADH}, \textit{Pyr}, \textit{Rib5P},
\textit{ADPf}, \textit{Phi})
is again an insight provided by the modularization.

The modularizations $\mathcal{T}_{\mathrm{MOD},1}$ and
$\mathcal{T}_{\mathrm{MOD},2}$ identified using the theory and methods developed
in the paper
constitute parsimonious, coarse-grained views of the metabolite
network studied and provide important insight concerning the dynamics
of the
biological system as a whole.

\section{Directions for future research}\label{sec7}

Application of the methods developed here to SKMs with large species
sets and many component reactions
is of considerable interest. In ongoing research that examines
biochemical signalling networks with approximately 900 reactions and
750 species, the methods have been found to work effectively and to
provide scientifically interesting modularizations.

It would be useful to consider methods for testing the adequacy of an SKM
(perhaps augmented to allow for measurement error) as a statistical
model of
a given cellular system. Testing conditional independence relationships
implied by a modularization of the SKM (such as the one in Figure \ref{TMod3}
for the red blood cell) offers a
promising means of assessing model adequacy. Clearly, it is not
necessary to
measure experimentally all species in the SKM, but all species in the
relevant conditioning set (separator) must be measured. Intuitively,
with $A\perp_{G^{\thicksim}}B|D$,
changes in $B$---perhaps resulting from direct intervention on the
levels of $B$---should be uninformative about changes in $A$ over time
intervals sufficiently short to ensure that levels of $D$ usually
remain constant (and vice versa).

SKMs subject to interventions are likely to become an area of active
research, given their relevance both to medical and biotechnological
applications. The predicted effect of interventions (e.g., gene
knock-outs, RNA silencing, or receptor inhibition)
could be derived by altering the specification of the SKM accordingly and
comparing with the original SKM. There are also interesting connections
with the causal inference
literature more generally. Recently, Commenges and G\'{e}gout-Petit \cite{Commenges2009} introduced a
``general dynamical
model as a framework for causal interpretation,'' adopting an approach to
causality based on ``physical laws in sufficiently large systems.'' Local
independence plays an important role in their
analysis and definition of influence. One might imagine that a sufficiently
large SKM would be a candidate ``perfect system'' for a given smaller and
observable cellular system. However, the jump processes followed by
biochemical species and hence also SKMs do not belong to the class ($%
\mathcal{D}$) of special semimartingales to which \cite{Commenges2009}
confines attention. Nevertheless, the approach seems relevant in broad
terms. Finally, the experimental design of interventions to test causal
claims derived from SKMs merits attention.

\begin{appendix}
\section{Proofs for global dynamic independence}\label{appA}

\subsection{\texorpdfstring{Proof of Theorem \protect\ref{Main1}}{Proof of Theorem 4.4}}\label{appA1}

First we show that $\mathcal{F}_{t}^{\Delta(A)\cap\Delta(B)}=\mathcal
{F}%
_{AB}^{D}(t)\subseteq\mathcal{F}_{t}^{D^{\ast}}$, in order to establish
(\ref{First}) below---that is, the internal history of all
reactions that change
$A$ and $B$ is contained in the internal history of $N^{D^{\ast}}(t)$.

The separation $A\bot_{G^{\thicksim}}B|D$ implies that for any $m\in
\Delta(A)\cap\Delta(B)$, $R[m]\subseteq D$ [suppose not---then in the
KIG $G$, either $\mathrm{pa}(B)\cap A\neq\varnothing$ or $\mathrm{pa}(A)\cap B\neq
\varnothing$
which contradicts the separation]. For any $m\in\Delta(A)\cap\Delta(B)$,
$R[m]\neq\varnothing$ and $R^{\ast}[m]\neq\varnothing$ by Definition
\ref{stdSKM}(iii) and (iv); clearly $R^{\ast}[m]\subseteq D$. Hence $m\in
\Delta(D)$ and $\Delta D_{AB}=\Delta(A)\cap\Delta(B)$ (with the
possibility $\Delta D_{AB}=\varnothing$ not excluded). By
Condition \ref{IdentGamma}, any reaction in $\Gamma= \Delta(A)\cap\Delta(B)$
changes $D$
differently---that is, the partition $\mathcal{M}(\Delta D_{AB})$ is either
empty or consists of singletons---since $\forall m,\tilde{m}\in\Gamma
$ $%
(m\neq\tilde{m})$, $S_{m}^{-}\neq S_{\tilde{m}}^{-}$ and $%
(S_{m}^{-})^{A}=(S_{m}^{-})^{B}=\mathbf{0}$ because $R^{\ast
}[m]\subseteq D$
(similarly for $\tilde{m}$), hence $(S_{m}^{-})^{D}\neq(S_{\tilde{m}%
}^{-})^{D}$ and $S_{m}^{D}\neq S_{\tilde{m}}^{D}$. Hence $N_{t}^{\Delta
(A)\cap\Delta(B)}=N_{AB}^{D}(t)$ $\forall t$ and $\mathcal
{F}_{t}^{\Delta
(A)\cap\Delta(B)}=\mathcal{F}_{AB}^{D}(t)\subseteq\mathcal{F}%
_{t}^{D^{\ast}}$, which implies immediately that
%
%
\begin{equation} \label{First}
\mathcal{F}_{t}^{\Delta(A)}\indep{\mathcal{F%
}_{t}^{\Delta(A)\cap\Delta(B)}}|\mathcal{F}_{t}^{D^{\ast}};\mathsf{%
\tilde{P}}_{t}.
\end{equation}

Together with
%
%
\begin{equation} \label{rem}
\mathcal{F}_{t}^{\Delta(A)}\indep{\mathcal{F%
}_{t}^{\Delta(B)\setminus\Delta(A)}}|\mathcal{F}_{t}^{D^{\ast
}};\mathsf{%
\tilde{P}}_{t}
\end{equation}
(which is proved below) it follows that
\[
\mathcal{F}_{t}^{\Delta(A)}\indep{\mathcal{F%
}_{t}^{\Delta(B)}}|\mathcal{F}_{t}^{D^{\ast}};\mathsf{\tilde{P}}_{t},
\]
since $\mathcal{F}_{t}^{\Delta(A)\cap\Delta(B)}\vee\mathcal{F}%
_{t}^{D^{\ast}}=\mathcal{F}_{t}^{D^{\ast}}$ and $\mathcal
{F}_{t}^{\Delta
(B)\setminus\Delta(A)}\vee\mathcal{F}_{t}^{\Delta(A)\cap\Delta
(B)}=%
\mathcal{F}_{t}^{\Delta(B)}$. It then follows that $\mathcal
{F}_{t}^{A}\indep\mathcal{F}_{t}^{B}|\mathcal{F}%
_{t}^{D^{\ast}};\mathsf{\tilde{P}}_{t}$ as required since it is clear from
the definition of $N^{A}(t)$ and $N^{B}(t) $ that $\mathcal{F}%
_{t}^{A}\subseteq\mathcal{F}_{t}^{\Delta(A)}$ and $\mathcal{F}%
_{t}^{B}\subseteq\mathcal{F}_{t}^{\Delta(B)}$. (The reader unfamiliar with
conditional independence of $\sigma$-fields and its properties is referred
to \cite{FMR90}---see, in the context of this proof, Theorem 2.2.1, Corollary
2.2.4, Theorem 2.2.10 and Corollary 2.2.11 there.)

It remains to establish (\ref{rem}). Under $\mathsf{\tilde{P}}$, and
hence also under $\mathsf{\tilde{P}}_{t}$, $\{\mathcal
{F}_{t}^{m}|m=1,\ldots,M
\} $ are independent $\sigma$-fields (see Lemma \ref{Likl}). It follows
that $\mathcal{F}_{t}^{\Delta(A)}\indep\break{%
\mathcal{F}_{t}^{\Delta(B)\setminus\Delta(A)}}|\mathcal
{F}_{t}^{\Delta
D_{D}};\mathsf{\tilde{P}}_{t}$ since $[\Delta(A),\Delta(B)\setminus
\Delta(A),\Delta D_{D}]$ is a partition of $\{1,\break \ldots,M\}$. It now suffices
for (\ref{rem}) to show the existence of $\mathcal{G}_{t}^{\Delta
(A)}\subseteq\mathcal{F}_{t}^{\Delta(A)}$ and $\mathcal{G}_{t}^{\Delta
(B)\setminus\Delta(A)}\subseteq\mathcal{F}_{t}^{\Delta(B)\setminus
\Delta(A)}$, such that $\mathcal{G}_{t}^{\Delta(A)}\vee\mathcal{G}%
_{t}^{\Delta(B)\setminus\Delta(A)}\vee\mathcal{F}_{t}^{\Delta
D_{D}}=%
\mathcal{F}_{t}^{D^{\ast}}$. Heuristically, we want to identify
``information'' contained only in $\mathcal{F}_{t}^{\Delta(A)}$ and $\mathcal{F}%
_{t}^{\Delta(B)\setminus\Delta(A)}$, respectively, which when jointly
combined with $\mathcal{F}_{t}^{\Delta D_{D}}$ gives the internal
history of
$N_{t}^{D^{\ast}}$. But this corresponds exactly to the way
$N_{t}^{D^{\ast
}}=[\{N_{A}^{D}(t),N_{AB}^{D}(t)\},N_{B}^{D}(t),N_{D}^{D}(t)]$ was constructed.

Recall Definition \ref{NDStar} for $N^{D^{\ast}}(t)$; its
history $\mathcal{F}_{t}^{D^{\ast}}$ is given by $[\mathcal{F}%
_{A}^{D}(t)\vee\mathcal{F}_{AB}^{D}(t)\vee\mathcal{F}_{B}^{D}(t)\vee
\mathcal{F}_{D}^{D}(t)]$. Since $\Delta D_{A}\subseteq\Delta(A)$, $%
\mathcal{F}_{A}^{D}(t)\subseteq\mathcal{F}_{t}^{\Delta D_{A}}\subseteq
\mathcal{F}_{t}^{\Delta(A)}$; similarly $\Delta D_{AB}\subseteq\Delta(A)$
and hence $\mathcal{F}_{AB}^{D}(t)=\mathcal{F}_{t}^{\Delta
D_{AB}}\subseteq
\mathcal{F}_{t}^{\Delta(A)}$; and $\Delta D_{B}\subseteq\Delta
(B)\setminus\Delta(A)$ hence $\mathcal{F}_{B}^{D}(t)\subseteq
\mathcal{F}%
_{t}^{\Delta D_{B}}\subseteq\mathcal{F}_{t}^{\Delta(B)\setminus
\Delta
(A)}$. Note that $\mathcal{F}_{D}^{D}(t)=\mathcal{F}_{t}^{\Delta D_{D}}$
since $\mathcal{M}(\Delta D_{D})$ is either empty or consists of
singletons---any 2 reactions that change $D$ alone must do so
differently since
no 2
columns of $S$ are equal (by Definition \ref{SKM}). Finally, taking $%
\mathcal{G}_{t}^{\Delta(A)}=\mathcal{F}_{A}^{D}(t)\vee\mathcal{F}%
_{AB}^{D}(t)\subseteq\mathcal{F}_{t}^{\Delta(A)}$ and $\mathcal{G}%
_{t}^{\Delta(B)\setminus\Delta(A)}=\mathcal{F}_{B}^{D}(t)\subseteq
\mathcal{F}_{t}^{\Delta(B)\setminus\Delta(A)}$ completes the proof since
then $\mathcal{G}_{t}^{\Delta(A)}\vee\mathcal{G}_{t}^{\Delta
(B)\setminus
\Delta(A)}\vee\mathcal{F}_{t}^{\Delta D_{D}}=\mathcal{F}_{t}^{D^{\ast}}$
as required.

\subsection{\texorpdfstring{Proof of Theorem \protect\ref{Main2}}{Proof of Theorem 4.5}}

The proof is in 3 parts.
(I) First show that \mbox{$\forall m\in\mathcal{M}$}, either $R[m]\subseteq
A\cup
D $ in which case $\int_{0}^{t}\lambda_{m}(u)\,du$ is adapted to
$\mathcal{F}%
_{t}^{AD^{\ast}}$, or $R[m]\subseteq B\cup D$ in which case $%
\int_{0}^{t}\lambda_{m}(u)\,du$ is adapted to $\mathcal{F}_{t}^{BD^{\ast}}$.

The separation $A\bot_{G^{\thicksim}}B|D$ implies that either $%
R[m]\subseteq A\cup D$ or $R[m]\subseteq B\cup D$. Suppose not, then
$B\cap
R[m]\neq\varnothing$ and $A\cap R[m]\neq\varnothing$---arguing using
(i) of
Definition \ref{stdSKM}, either $m\in\Delta(A)$ in which case $B\cap
\mathrm{pa}(A)\neq\varnothing$, which contradicts the separation$;$ or $m\in
\Delta
(B)$ in which case $A\cap \mathrm{pa}(B)\neq\varnothing$, which also contradicts the
separation. If $B\cap R[m]\neq\varnothing$ and $A\cap R[m]\neq
\varnothing$,
then $m\in\Delta D_{D}$ is not possible---if $m\in\Delta D_{D}$ then the
reactants that are changed $R^{\ast}[m]\subseteq D$ and hence, by (iv) of
Definition \ref{stdSKM}, either $B\cap R[m]\neq\varnothing$ or$A\cap
R[m]\neq\varnothing$ but not both.

Therefore, if $R[m]\subseteq A\cup D$ (resp., $R[m]\subseteq B\cup D$) then
both $\lambda_{m}(t)$ and $\log(\lambda_{m}(t))$ are measurable with
respect to $\mathcal{F}_{t}^{R[m]}\subseteq\mathcal{F}_{t}^{A\cup
D}\subseteq$ $\mathcal{F}_{t}^{AD^{\ast}}$ (resp., $\mathcal
{F}_{t}^{B\cup
D}\subseteq\mathcal{F}_{t}^{BD^{\ast}}$) by (\ref{cmgm}),
since $%
X^{R[m]}(t-)$ is measurable $\mathcal{F}_{t}^{R[m]}$ and $\mathcal{F}%
_{t}^{D}\subseteq\mathcal{F}_{t}^{D^{\ast}}$. Since $\lambda
_{m}(t)$ is also c\`{a}gl\`{a}d, if $R[m]\subseteq A\cup D$ (resp., $%
R[m]\subseteq B\cup D$) then $\lambda_{m}(t)$ is $\mathcal
{F}_{t}^{AD^{\ast
}}$-predictable and hence $\int_{0}^{t}\lambda_{m}(u)\,du$ is $\mathcal
{F}%
_{t}^{AD^{\ast}}$-adapted (resp., $\mathcal{F}_{t}^{BD^{\ast}}$-adapted).

(II) Second show that if $R[m]\subseteq A\cup D$ (resp., $R[m]\subseteq
B\cup
D$), then $\{T_{s}^{m}\}_{s\geq1}$ are $\mathcal{F}_{t}^{AD^{\ast}}$%
-stopping times (resp., $\mathcal{F}_{t}^{BD^{\ast}}$-stopping times). It
will then follow that $1(T_{s}^{m}\leq t)\log(\lambda
_{m}(T_{s}^{m}))$ is $%
\mathcal{F}_{t}^{AD^{\ast}}$-measurable $\forall s\geq1$ (resp.,
$\mathcal{F%
}_{t}^{BD^{\ast}}$-measurable) by the definition of $\mathcal{F}%
_{T_{s}^{m}}^{AD^{\ast}}$, because $\log(\lambda_{m}(t))$ is
left continuous and hence $\log(\lambda_{m}(T_{s}^{m}))$ is $\mathcal
{F}%
_{T_{s}^{m}}^{AD^{\ast}}$-measurable---see, for example, Theorem 2.1.10 of
\cite{LastBrandt95}.
This in turn yields that $\sum_{s\geq1}1(T_{s}^{m}\leq
t)\log(\lambda_{m}(T_{s}^{m}))$ is $\mathcal{F}_{t}^{AD^{\ast}}$%
-measurable (resp., $\mathcal{F}_{t}^{BD^{\ast}}$-measurable). To establish
the required stopping time property for $\{T_{s}^{m}\}_{s\geq1}$,
distinguish the following cases, exactly one of which must hold $\forall
m\in\mathcal{M}$:

(i) $m\in\Delta D_{D}$: recall that $\mathcal{M}(\Delta D_{D})$
consists of
singletons since any 2 reactions that change $D$ alone must do so
differently (by Definition \ref{SKM}). Therefore $N_{m}(t)$ is adapted
to $%
\mathcal{F}_{t}^{\Delta D_{D}}=\mathcal{F}_{D}^{D}(t)\subseteq\mathcal
{F}%
_{t}^{D^{\ast}}$, hence $\{T_{s}^{m}\}_{s\geq1}$ are $\mathcal{F}%
_{t}^{D^{\ast}}$-stopping times. Either $R[m]\subseteq A\cup D$ or $%
R[m]\subseteq B\cup D$ [by part (I) above]. If $R[m]\subseteq A\cup D$
(resp., $R[m]\subseteq B\cup D$), then $\{T_{s}^{m}\}_{s\geq1}$ are
necessarily $%
\mathcal{F}_{t}^{AD^{\ast}}$-stopping times (resp., $\mathcal{F}%
_{t}^{BD^{\ast}}$-stopping times), as required.

\mbox{\phantom{i}}(ii) $m\in\Delta A\cap\Delta(B)$: recall from the proof of Theorem
\ref{Main1} that $R[m]\subseteq D$ and the partition $\mathcal{M}(\Delta D_{AB})$
consists of singletons. Hence $\mathcal{F}_{t}^{m}\subseteq\mathcal{F}%
_{t}^{\Delta A\cap\Delta(B)}=\mathcal{F}_{AB}^{D}(t)\subseteq\mathcal
{F}%
_{t}^{D^{\ast}}$, $N_{m}(t)$ is adapted to $\mathcal{F}_{t}^{D^{\ast}}$
and $\{T_{s}^{m}\}_{s\geq1}$ are $\mathcal{F}_{t}^{D^{\ast}}$-stopping
times.

(iii) $m\in\Delta(A)\setminus\Delta(B)$: we have that $R[m]\subseteq
A\cup D$ by part (I) above; consider the cases (iiia) $m\in\Delta
D_{A}$, and
(iiib) $m\notin\Delta D_{A}$ in turn below to conclude that in each
case $%
\{T_{s}^{m}\}_{s\geq1}$ are $\mathcal{F}_{t}^{AD^{\ast}}$-stopping times.

\mbox{\phantom{(iii)} }(iiia) Identify the element of $\mathcal{M}(\Delta D_{A})$ corresponding to
changes in $D$ equal to $S_{m}^{D}$, $\mathcal{M}_{e}(\Delta D_{A})$ say.
Denote the corresponding element of $N_{A}^{D}(t)$ by $N_{A,e}^{D}(t)$,
which is clearly measurable $\mathcal{F}_{t}^{D^{\ast}}$, and the jump
times of this univariate counting process by $\{T_{A,e}^{D}(s)\}_{s\geq
1}$. Thus $%
\{T_{A,e}^{D}(s)\}_{s\geq1}$ are $\mathcal{F}_{t}^{AD^{\ast}}$-stopping
times. Now $\mathcal{M}_{e}(\Delta D_{A})$ may not be a singleton, but we
can write
\[
N_{m}(t)=\sum_{s\geq1}1\{T_{A,e}^{D}(s)\leq t\}1\{X^{A\cup
D}(T_{A,e}^{D}(s))-X^{A\cup D}(T_{A,e}^{D}(s)-)=S_{m}^{A\cup D}\},
\]
since $\nexists\tilde{m}\in\mathcal{M}_{e}(\Delta D_{A})$ $(m\neq
\tilde{m}%
)$ s.t. $S_{m}^{A\cup D}=S_{\tilde{m}}^{A\cup D}$ (by Definition \ref{SKM}
and $S_{m}^{B}=S_{\tilde{m}}^{B}=\mathbf{0}$). Since $X^{A\cup D}(t)$
is right
continuous and $X^{A\cup D}(t-)$ left continuous, and both are $\mathcal
{F}%
_{t}^{AD^{\ast}}$-adapted since $\mathcal{F}_{t}^{A\cup D}$-adapted, $%
[X^{A\cup D}(T_{A,e}^{D}(s))-X^{A\cup D}(T_{A,e}^{D}(s)-)]$ is $\mathcal
{F}%
^{AD^{\ast}}(T_{A,e}^{D}(s))$-measurable (by, e.g., Theorem 2.1.10 of
\cite{LastBrandt95}). Hence the summand is $\mathcal{F}^{AD^{\ast
}}(t)$-measurable $\forall s\geq1$, $N_{m}(t)$ is adapted to $\mathcal
{F}%
_{t}^{AD^{\ast}}$ and $\{T_{s}^{m}\}_{s\geq1}$ are $\mathcal{F}%
_{t}^{AD^{\ast}}$-stopping times.

\mbox{\phantom{(iii)} }(iiib) Then $m\in\Delta^{\ast}(A):=\Delta(A)\setminus(\Delta
(B)\cup
\Delta(D))$. Define for any $p$-variate counting process $N(t)$
$(p\geq1)$, the ``ground
process'' $\bar{N}(t):=\mathbf{1}_{p\times1}^{\prime}N(t)$. We may then
write
\[
\bar{N}^{\Delta^{\ast}(A)}(t)=\bar{N}^{A}(t)-\bar{N}_{A}^{D}(t)-\bar
{N}%
_{AB}^{D}(t),
\]
where $\bar{N}^{\Delta^{\ast}(A)}(t)$ is the number of reactions on $[0,t]$
that change $A$ alone [noting that $\Delta(A)\cap\Delta(B)=\Delta
D_{AB}$]. Hence $\bar{N}^{\Delta^{\ast}(A)}(t)$ is measurable $\mathcal{F}%
_{t}^{A}\vee\mathcal{F}_{A}^{D}(t)\vee\mathcal{F}_{AB}^{D}(t)\subseteq
\mathcal{F}_{t}^{AD^{\ast}}$, and its jump times $\{T_{s}^{\Delta
^{\ast
}(A)})\}_{s\geq1}$ are $\mathcal{F}_{t}^{AD^{\ast}}$-stopping times. Also,
\[
N_{m}(t)=\sum_{s\geq1}1\bigl\{T_{s}^{\Delta^{\ast}(A)}\leq
t\bigr\}1\bigl\{X^{A}\bigl(T_{s}^{\Delta^{\ast}(A)}\bigr)-X^{A}\bigl(T_{s}^{\Delta^{\ast
}(A)}-\bigr)=S_{m}^{A}\bigr\},
\]
since $\nexists\tilde{m}\in\Delta^{\ast}(A)$ $(m\neq\tilde{m})$
s.t. $%
S_{m}^{A}=S_{\tilde{m}}^{A}$ (by Definition \ref{SKM} and $S_{m}^{B\cup
D}=S_{\tilde{m}}^{B\cup D}=\mathbf{0}$). Since $X^{A}(t)$ is right continuous
and $X^{A}(t-)$ left continuous, and both are $\mathcal
{F}_{t}^{AD^{\ast}}$%
-adapted since $\mathcal{F}_{t}^{A}$-adapted, $[X^{A}(T_{s}^{\Delta
^{\ast
}(A)})-X^{A}(T_{s}^{\Delta^{\ast}(A)}-)]$ is $\mathcal{F}^{AD^{\ast
}}(T_{s}^{\Delta^{\ast}(A)})$-measurable. Hence the summand is
$\mathcal{F}%
^{AD^{\ast}}(t)$-measurable $\forall s\geq1$, and $\{T_{s}^{m}\}
_{s\geq1}$
are $\mathcal{F}_{t}^{AD^{\ast}}$-stopping times.

(iv) $m\in\Delta(B)\setminus\Delta(A)$: we have that $R[m]\subseteq
B\cup D$ by part (I) above; argue as in (iii) with $A$ in place of $B$ and
vice versa to conclude that $\{T_{s}^{m}\}_{s\geq1}$ are $\mathcal{F}%
_{t}^{BD^{\ast}}$-stopping times.

(III) Combining parts (I) and (II) above establishes that if
$R[m]\subseteq
A\cup D$ (resp., $R[m]\subseteq B\cup D$) then $\mathcal{L}%
_{m,t}:=\exp(l_{m}(t))$ is measurable $\mathcal{F}_{t}^{AD^{\ast}}$
(resp., $%
\mathcal{F}_{t}^{BD^{\ast}}$). Then, in an obvious manner, grouping
the $%
\mathcal{L}_{m,t}$ into 2 groups according to the forementioned
measurability property and defining the $\psi_{iD^{\ast},t}$ as the
product within each group yields $\mathcal{L}_{t}=\psi_{AD^{\ast
},t}\cdot
\psi_{BD^{\ast},t}$, where $\psi_{iD^{\ast},t}$ is nonnegative and
$%
\mathcal{F}_{t}^{i}\vee\mathcal{F}_{t}^{D^{\ast}}$-measurable for
$i\in
\{A,B\}$.

\section{Additional proofs}\label{Prooffs}

\vspace*{-10pt}

\begin{pf*}{Proof of Lemma \ref{Filts}}
It remains to establish that
$\mathcal{F}_{t}^{A}=\mathcal{F}_{t}^{X^{A}}$. First show that $\mathcal
{F}_{t}^{A}\supseteq\mathcal{F}_{t}^{X^{A}}$. We have that $\mathcal
{F}_{t}^{A}=\sigma(Z_{s}^{A}1(T_{s}^{A}\leq u);0\leq
u\leq t$, $s\geq1)$ and $X^{A}(u)=X^{A}(0)+\sum_{s\geq
1}Z_{s}^{A}1(T_{s}^{A}\leq u)$, which is therefore measurable $\mathcal
{F}%
_{t}^{A}$. Second show that $\mathcal{F}_{t}^{A}\subseteq\mathcal{F}%
_{t}^{X^{A}}$. We have also that $\mathcal{F}_{t}^{A}=\sigma
(1(Z_{s}^{A}=S_{m}^{A})1(T_{s}^{A}\leq u);0\leq u\leq t$, $s\geq1,m\in
\Delta A)$, hence it suffices to show that $%
1(Z_{s}^{A}=S_{m}^{A})1(T_{s}^{A}\leq u)$ is measurable $\mathcal{F}%
_{t}^{X_{A}}$. By its construction, $\{T_{s}^{A}\}$ are the jump times of
the right-continuous jump process $X^{A}$. The filtration $\{\mathcal{F}
_{t}^{X_{A}}\}$ is right continuous. Hence, for $s\geq1$, $T_{s}^{A}$
is an
$\mathcal{F}_{t}^{X_{A}}$-stopping time and $X^{A}(T_{s}^{A})$ is
$\mathcal{F%
}^{X_{A}}(T_{s}^{A})$-measurable. Since $%
Z_{s}^{A}=X^{A}(T_{s}^{A})-X^{A}(T_{s-1}^{A})$ and $\mathcal{F}%
^{X_{A}}(T_{s-1}^{A})\subseteq\mathcal{F}^{X_{A}}(T_{s}^{A})$, $Z_{s}^{A}$
is also $\mathcal{F}^{X_{A}}(T_{s}^{A})$-measurable. Hence $%
\{1(Z_{s}^{A}=S_{m}^{A})=1\}\cap\{1(T_{s}^{A}\leq u)=1\}\in\mathcal{F}
^{X_{A}}(u)\subseteq\mathcal{F}^{X_{A}}(t)$ by the definition of
$\mathcal{F%
}^{X_{A}}(T_{s}^{A})$, and therefore
$1(Z_{s}^{A}=S_{m}^{A})1(T_{s}^{A}\leq
u)$ is measurable $\mathcal{F}_{t}^{X_{A}}$.
\end{pf*}
\begin{pf*}{Proof of Lemma \ref{CIDom}}
Let $\mathcal{L}_{i3}:=(\mathsf
{dP}/\mathsf{d\tilde{P}})|_{\mathcal{F}%
^{i}\vee\mathcal{F}^{3}}$, and $\mathcal{L}_{3}:=(\mathsf{dP}/\mathsf{d
\tilde{P}})|_{\mathcal{F}^{3}}$. Then it is straightforward to show
that $%
\mathcal{L}_{i3}=\mathsf{\tilde{E}}[\mathcal{L}_{123}|\mathcal
{F}^{i}\vee
\mathcal{F}^{3}]$ and $\mathcal{L}_{3}=\mathsf{\tilde{E}}[\mathcal
{L}_{123}|%
\mathcal{F}^{3}]$, where $\mathsf{\tilde{E}}$ denotes expectation under
$%
\mathsf{\tilde{P}}$. Hence, $\mathcal{L}_{13}=\psi_{13}\mathsf{\tilde
{E}}%
[\psi_{23}|\mathcal{F}^{1}\vee\mathcal{F}^{3}]$ and $\mathcal
{L}_{23}=\psi
_{23}\mathsf{\tilde{E}}[\psi_{13}|\mathcal{F}^{1}\vee\mathcal
{F}^{3}]$ by
the nonnegativity and measurability of the $\psi_{i3}$. Since
$\mathcal{F}%
^{2}\vee\mathcal{F}^{3}\indep{\mathcal{F}^{1}}|\mathcal{F}^{3};\mathsf{
\tilde{P}}$, $\mathsf{\tilde{E}}[\psi_{23}|\mathcal{F}^{1}\vee\mathcal
{F}%
^{3}]=\mathsf{\tilde{E}}[\psi_{23}|\mathcal{F}^{3}]$ by Definition \ref{CI}
and hence $\mathcal{L}_{13}=\psi_{13}\mathsf{\tilde{E}}[\psi
_{23}|\mathcal{%
F}^{3}]$. Similarly, $\mathcal{L}_{23}=\psi_{23}\mathsf{\tilde{E}}[\psi
_{13}|\mathcal{F}^{3}]$. Furthermore, $\mathcal{L}_{3}=\mathsf{\tilde
{E}}%
[\psi_{13}|\mathcal{F}^{3}]\mathsf{\tilde{E}}[\psi_{23}|\mathcal{F}^{3}]$
by the nonnegativity and measurability of the $\psi_{i3}$ and since $%
\mathcal{F}^{1}\vee\mathcal{F}^{3}\indep{\mathcal{F}^{2}\vee\mathcal
{F}%
^{3}}|\mathcal{F}^{3};\mathsf{\tilde{P}}$. Therefore,
%
%
\begin{equation}
\mathcal{L}_{123}\mathcal{L}_{3}=\mathcal{L}_{13}\mathcal{L}_{23}
\end{equation}
and, in particular, $\mathcal{L}_{123}\mathcal{L}_{3}=\mathcal{L}_{13}%
\mathcal{L}_{23}$ on the event $\{\mathcal{L}_{3}=\mathsf{\tilde
{E}}[\psi
_{13}|\mathcal{F}^{3}]\mathsf{\tilde{E}}[\psi_{23}|\mathcal{F}^{3}]>0\}$,
whence $\mathcal{F}^{1}\indep{\mathcal{F}^{2}}|\mathcal{F}^{3};\mathsf{P}$
by Theorem 2.2.14 of \cite{FMR90}.
\end{pf*}
\begin{pf*}{Proof of Proposition \ref{JTProp}} The proof is in 3 steps,
according to the number of pairs of clusters aggregated under step 5 of
Algorithm \ref{Tmod}: (i)
for the case where no pair of clusters is aggregated, and hence $%
\mathcal{T}_{\mathrm{MOD}}=\mathcal{T}_{C};$ (ii) for the case where exactly 1 pair
of clusters is aggregated; (iii) for the case where more than 1 pair of
clusters is aggregated.

{\smallskipamount=0pt
\begin{longlist}
\item
$\mathcal{T}_{C}$ is a junction tree representation of the
clique decomposition of $G_{T}^{\thicksim}$. For the proof of
this case see the proof of Theorem 4.6 of \cite{CDLS2007}.

\item
$\mathcal{T}_{\mathrm{MOD}}$ is connected (as a consequence of $\mathcal{%
T}_{C}$ being connected), and has $(\delta-1)$ nodes and $(\delta-2)$
edges (one less edge than $\mathcal{T}_{C}$); $\mathcal{T}_{\mathrm{MOD}}$ is
therefore a tree, whence there is a unique path in $\mathcal{T}_{\mathrm{MOD}}$
between any pair $(M_{d},M_{e})$ of its clusters. It is
straightforward (but somewhat tedious) to show that every cluster on this
path must contain $M_{d}\cap M_{e}$ since the corresponding path in
$\mathcal{T}_{C}$ possesses this junction property [by (i) above].
Hence $%
\mathcal{T}_{\mathrm{MOD}}$ is a junction tree. It remains to prove that for any 2
adjacent clusters $(M_{d},M_{e})$ in $\mathcal{T}_{\mathrm{MOD}}$, we have $%
M_{d}\cap M_{e}=V_{\mathit{de}}\cap V_{\mathit{ed}}$ and $V_{\mathit{de}}\perp_{G_{T}^{\thicksim%
}}V_{\mathit{ed}}|S_{\mathit{de}}$.

We will show (iia) that edges ``in common'' between $\mathcal{T}_{C}$
and $\mathcal{T}_{\mathrm{MOD}}$---the $(\delta-2)$ edges not
removed by the cluster aggregation---carry the same label, that is, the
intersection of the clusters joined by each such edge is unchanged; and
(iib) that cutting any such edge in both $\mathcal{T}_{C}$ and $\mathcal
{T}_{\mathrm{MOD}}$
results in pairs of subtrees whose clusters have identical unions in
the two cases.
The result then follows from (i) above.

\begin{longlist}[(iia)(ii)\mbox{\phantom{0(}}]
\item[(iia)] If both clusters, $(M_{d},M_{e})$, joined by such an edge, are in
$\mathcal{T}_{C}$ and $\mathcal{T}_{\mathrm{MOD}}$ the claim is obviously
true. Consider then the case where $M_{d}$, say, is the result of the
aggregation of the cluster pair $(M_{\alpha},M_{\beta})$. Suppose,
without loss of generality, that $M_{\alpha}\sim M_{e}$ in
$\mathcal{T}_{C}$. Now $S_{\mathit{de}}=(M_{e}\cap M_{\alpha})\cup
(M_{e}\cap M_{\beta})$. The edge joining $M_{d}$ to $M_{e}$ in
$\mathcal{T}_{\mathrm{MOD}}$ was formerly, in $\mathcal{T}_{C}$, the edge
$M_{\alpha
}\sim M_{e}$, whence $S_{\mathit{de}}=(M_{e}\cap M_{\alpha})$ since $M_{\alpha
}$ is
on the path between $M_{\beta}$ and $M_{e}$ in $\mathcal{T}_{C}$
and $(M_{e}\cap M_{\beta})\subseteq M_{\alpha}$. Thus, the
intersection of the
clusters joined by the edge is always the same in $\mathcal{T}_{\mathrm{MOD}}$
and $\mathcal{T}_{C}$, as claimed.

\item[(iib)] Let the edge that is cut in both cases be $M_{d}\sim M_{e}$
[where it is understood
that $M_{d}$, say, may be equal to $M_{\alpha}$ in $\mathcal{T}_{C}$
and hence equal to $(M_{\alpha}\cup M_{\beta})$ in $\mathcal{T}%
_{\mathrm{MOD}}$]. It is required to show, using an obvious notation, that $%
V_{\mathit{de}}^{\mathrm{MOD}}=V_{\mathit{de}}^{C}$ and $V_{\mathit{ed}}^{\mathrm{MOD}}=V_{\mathit{ed}}^{C}$. It is well known
that cutting an edge in any tree results in 2 disconnected subtrees.
One of
the 2 pairs of subtrees generated here must contain 2 identical
subtrees. Suppose then,
without loss of generality, that $\mathcal{T}_{\mathrm{MOD}}^{\mathit{ed}}=\mathcal{T}%
_{C}^{\mathit{ed}}$, whence $V_{\mathit{ed}}^{\mathrm{MOD}}=V_{\mathit{ed}}^{C}$. The subtrees $%
\mathcal{T}_{\mathrm{MOD}}^{\mathit{de}}$ and $\mathcal{T}_{C}^{\mathit{de}}$ have the same
clusters, except
for the aggregation of the cluster pair $(M_{\alpha},M_{\beta})$ to form
$M_{\alpha\beta}$ in $\mathcal{T}_{\mathrm{MOD}}^{\mathit{de}}$. It is
straightforward (but tedious) to show that, for $\gamma\notin\{\alpha
,\beta\}$ and $M_{\gamma}$ a cluster in $\mathcal{T}_{C}^{\mathit{de}}$,
\[
\lbrack M_{\gamma}]_{\mathcal{T}_{C}^{\mathit{de}}}\setminus\{M_{\alpha
},M_{\beta
}\}=[M_{\gamma}]_{\mathcal{T}_{\mathrm{MOD}}^{\mathit{de}}}\setminus\{M_{\alpha\beta
}\},
\]
where $[M]_{\mathcal{T}}$ are the clusters that can be reached from
cluster $%
M$ by paths in a tree~$\mathcal{T}$. The subtrees $\mathcal{T}_{\bullet
}^{\mathit{de}}$ are themselves connected graphs, but disconnected from the
corresponding $\mathcal{T}_{\bullet}^{\mathit{ed}}$. Therefore the clusters of $
\mathcal{T}_{\bullet}^{\mathit{de}}$ are given exactly by $[M_{\gamma
}]_{\mathcal{T}%
_{\bullet}^{\mathit{de}}}$, where $M_{\gamma}$ is any one of its clusters. It
follows that
\begin{eqnarray*}
V_{\mathit{de}}^{C}&=&\Bigl\{\bigcup[M_{\gamma}]_{\mathcal{T}_{C}^{\mathit{de}}}\setminus
\{M_{\alpha},M_{\beta}\}\Bigr\}\cup M_{\alpha}\cup M_{\beta}
\\
&=&\Bigl\{\bigcup[M_{\gamma}]_{\mathcal{T}_{\mathrm{MOD}}^{\mathit{de}}}\setminus\{M_{\alpha
\beta
}\}\Bigr\}\cup M_{\alpha\beta}=V_{\mathit{de}}^{\mathrm{MOD}}
\end{eqnarray*}
as required.
\end{longlist}

\item[(iii)]
The proof is by induction on the number of cluster pairs, $n$ say,
that are aggregated. Parts (i) and (ii) above establish the proposition
for $%
n=0$ and $n=1$. Exactly the same mode of argument as the one used in
(ii) above
also establishes that if the proposition holds for $n\geq0$, it must
hold for
$(n+1)$. This completes the proof.\qed
\end{longlist}}
\noqed\end{pf*}

\vspace*{-10pt}

\section{Species names for red blood cell SKM}
\label{species}

\textit{AMPf}${}={}$AMP (unbound); \textit{ADPf}${}={}$ADP;
\textit{ATPf}${}={}$ATP; \textit{DHAP}${}=$\break Dihydroxyacetone
phosphate; \textit{E4P}${}={}$Erythrose
4-phosphate;
\textit{Fru6P}${}={}\break$Fructose 6-phosphate; \textit{Fru16P2}${}={}$Fructose
1,6-phosphate; \textit{GlcA6P}${}={}\break$Phospho-D-glucono-1,5-lactone;
\textit{Glcin}${}={}$Glucose
(cytoplasmic); \textit{Glcout}${}={}$External Glucose; \textit
{Glc6P}${}={}$Glucose 6-phosphate; \textit{GraP}${}={}$Glyceraldehyde 3-phosphate;
\textit{Gri13P2}${}={}$1,3-Bisphospho-D-glycerate;
\textit{Gri3P}${}={}$3-Phospho-D-glycerate; \textit{Gri23P2}${}={}$2,3-Bisphospho-D-glycerate;
\textit{Gri2P}${}={}$2-Phospho-D-glycerate; \textit{GSH}${}={}$Reduced Glutathione;
\textit{GSSG}${}={}$Oxidized Glutathione; \textit{Lac}${}={}$Lactate;
\textit{Lacex}${}=$ External Lactate; \textit{MgATP}; \textit{MgADP};
\textit{MgAMP}; \textit{Mg}; \textit{MgGri23P2}; \textit{NADH};
\textit{NADPf}${}={}$NADP (unbound); \textit{NADPHf}${}={}$NADPH;
\textit{P1f}${}={}$Protein1; \textit{P2}${}={}$Protein2;
\textit{P1NADP}${}={}$Protein1 bound NADP;\break
\textit{P1NADPH}${}={}$Protein1 bound NADPH; \textit{P2NADP}${}={}$Protein2
bound NADP;
\textit{P2NADPH}${}={}$Protein2 bound NADPH; \textit{PEP}${}={}$Phosphoenolpyruvate;
\textit{Phi}${}={}$Phosphate; \textit{NAD};
\textit{PRPP}${}=$ Phosphoribosylpyrophosphate; \textit{Pyr}${}={}$Pyruvate;
\textit{Pyrex}${}={}$External Pyruvate; \textit{Rib5P}${}={}$Ribose
5-phosphate; \textit{Rul5P}${}={}$Ribulose 5-phosphate;
\textit{Sed7P}${}=$ Sedoheptulose 7-phosphate;
\textit{Xul5P}${}={}$Xylulose\break 5-phosphate.
\end{appendix}

\section*{Acknowledgments}
The author is grateful for the research environment provided by the
Statistical Laboratory and the Cambridge Statistics Initiative, and to
A. P. Dawid, V. Didelez, C. Holmes, M. Jacobsen and D. J. Wilkinson for
helpful discussions. The comments of the Associate Editor and referees
were valuable in improving the paper. Computations were performed using
R version 2.8.1 and R packages gRbase and Rgraphviz.

\printaddresses

\end{document}